\documentclass[preprint,12pt]{elsarticle}

\usepackage{amssymb}
\usepackage{tabularx} 
\usepackage{booktabs}
\usepackage{amsmath} 
\usepackage{xcolor}
\usepackage{url}
\usepackage[section]{placeins}
\usepackage{xcolor}

\journal{International Journal of Electrical Power \& Energy Systems}

\begin{document}

\begin{frontmatter}





\title{Debunking the Speed-Fidelity Trade-Off: Speeding-up Large-Scale Energy Models while Keeping Fidelity}


\author[inst1,inst4]{Diego A. Tejada-Arango}
\author[inst3,inst5]{Juha Kiviluoma}
\author[inst1,inst2]{Germán Morales-España}

\affiliation[inst1]{organization={Energy \& Materials Transition, Netherlands Organisation for Applied Scientific Research (TNO)},
            country={The Netherlands}}

\affiliation[inst2]{organization={Faculty of Electrical Engineering, Mathematics and Computer Science, Delft University of Technology (TUDelft)},
            country={The Netherlands}}

\affiliation[inst3]{organization={Design and Operation of Energy Systems, VTT Technical Research Centre of Finland Ltd.},
            country={Finland}}

\affiliation[inst4]{organization={Instituto de Investigación Tecnológica, Escuela Técnica Superior de Ingeniería, Universidad Pontificia Comillas},
            country={Spain}}

\affiliation[inst5]{organization={Nodal-Tools Ltd.},
            country={Finland}}

\begin{abstract}

Energy system models are essential for planning and supporting the energy transition. However, increasing temporal, spatial, and sectoral resolutions have led to large-scale linear programming (LP) models that are often (over)simplified to remain computationally tractable—frequently at the expense of model fidelity. This paper challenges the common belief that LP formulations cannot be improved without sacrificing their accuracy. 
Inspired by graph theory, we propose to model energy systems using \textit{energy assets} (vertices), as a single building-block, and \textit{flows}  to connect between them. This reduces the need for additional components such as nodes and connections. The resulting formulation is more compact, without sacrificing accuracy, and leverages the inherent graph structure of energy systems.
To evaluate performance, we implemented and compared four common modelling approaches varying in their use of building blocks and flow representations. We conducted experiments using TulipaEnergyModel.jl and applied them to a multi-sector case study with varying problem sizes.
Results show that our single-building-block (1BB-1F) approach reduces variables and constraints by 26\% and 35\%, respectively, and achieves a 1.27x average speedup in solving time without any loss in model fidelity. The speedup increases with problem size, making this approach particularly advantageous for large-scale models.
Our findings demonstrate that not all LPs are equal in quality and that better reformulations can lead to substantial computational benefits. This paper also aims to raise awareness of model quality considerations in energy system optimisation and promote more efficient formulations without compromising fidelity.

\end{abstract}

\begin{graphicalabstract}
\includegraphics[width=\textwidth]{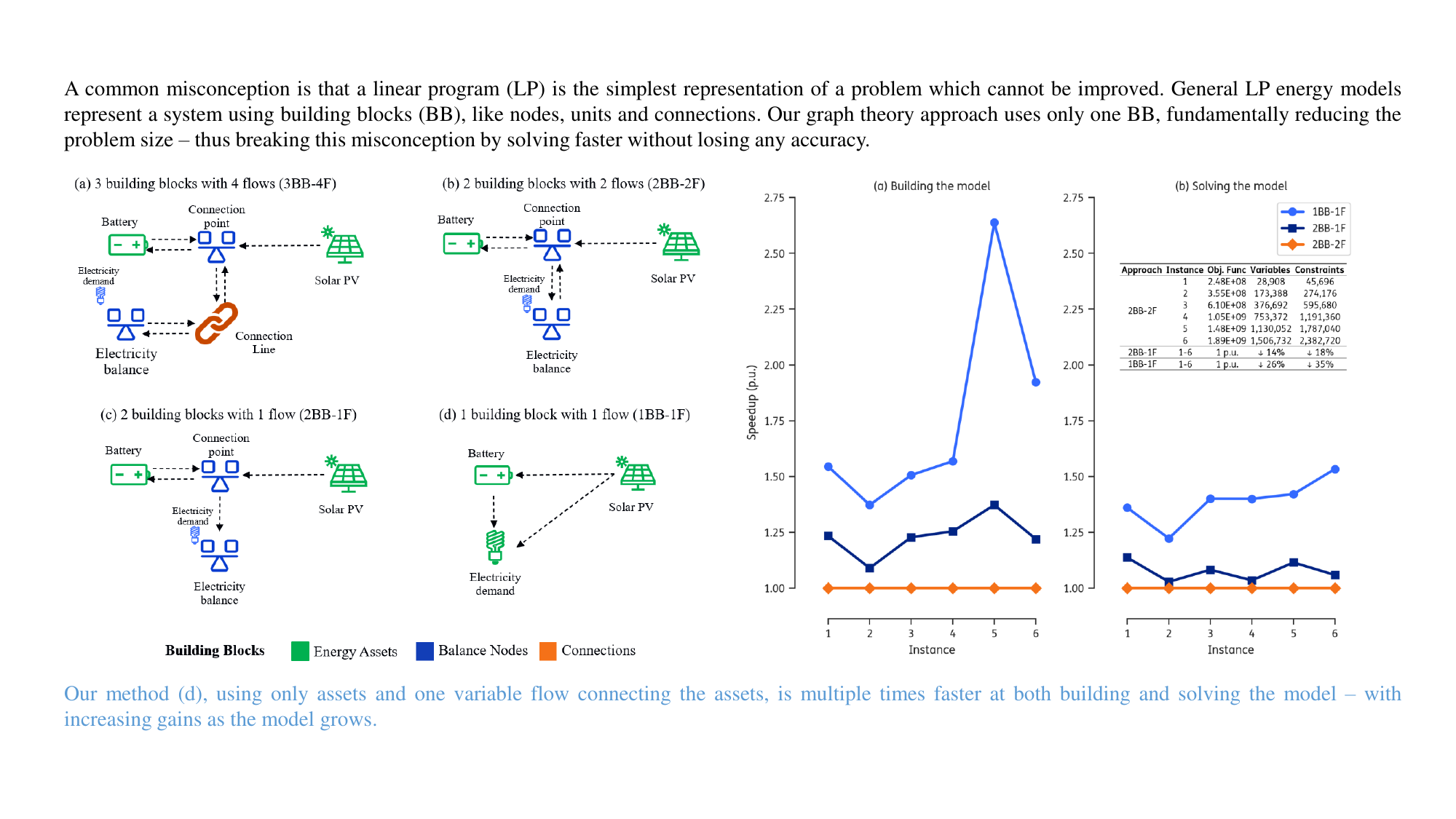}
\end{graphicalabstract}


\begin{highlights}

\item The proposed single modelling block reduces size, speeding up creation and runtime
\item Fewer variables and constraints without losing accuracy in energy optimisation models
\item Benefits increase with problem size, aiding large-scale energy system modelling
\item Enhancing modelling with a more direct and flexible connection between assets

\end{highlights}


\end{frontmatter}



\section{Introduction}
\label{sec:introduction}

\FloatBarrier
\subsection{Motivation}
\label{sec:motivation}
Although currently available models can integrate the energy sector coupling in the models in different ways, the implications of having flexibility in the modelling choices have not been given much attention. Computational performance is affected by these modelling choices, especially in large-scale energy problems like those in regional-level studies (e.g., European case studies) for long-term expansion planning (e.g., pathway planning). As the energy modelling community is increasingly exploring new system configurations to integrate different processes, such as green hydrogen production, hybrid operation between a storage asset and a renewable asset, or small modular reactors producing both electricity and heat, there is also a need for general and flexible modelling structures that can be adapted to new configurations without compromising the model's performance. Therefore, it is important to identify the computational implications of different modelling choices in this new context.

\subsection{Solving methods in energy systems modelling}
Many researchers view Linear Programming (LP) as the simplest problem representation that cannot be improved without sacrificing accuracy in energy systems modelling. This is a common misconception, which can lead to the creation of large overhead model sizes. There are at least two alternatives to improve the representation without losing accuracy. First, LP models can be reformulated to leverage specific characteristics, such as sparsity, or tailored for specialised algorithms, thus enhancing computational efficiency without sacrificing accuracy. For example, specialised algorithms like the Hungarian Algorithm for the assignment problem \cite{Kuhn2010} and the Push-Relabel method for network flow problems \cite{Cherkassky1997} often perform better than general LP solvers. These cases suggest that researchers should explore specialised methods to improve efficiency rather than relying solely on advancements in LP solvers or hardware. Second, since researchers widely use LP models, we explore another alternative to improve the formulation that still uses LP: reformulating the problem to exploit the graph nature of energy system models. All models are not the same, and improving the quality of an LP model means that it retains fidelity while solving faster. Hence, proposing a reformulation that can result in more significant advantages as the problem size increases is especially important for large-scale energy problems.

\subsection{Flexible modelling in energy system models}
\label{sec:background}
Multi-physics energy system models can be formulated in many ways, ranging from process-specific equations to a more generic approach where concepts such as nodes and units represent a wide variety of conversion, production, consumption, transfer, and storage processes \citep{ESDL2020}. In the process-specific approach, each process is described by its equation or set of equations. Models may have been built this way for historical reasons: at first the intention was to model a specific sector, such as the power system, and only afterwards was the model expanded to other energy sectors. An example is the early open-source model Balmorel \citep{Wiese2018}, which was initially a power system model but has since been expanded to most energy sectors. This trend has also been seen in commercial options like Plexos \citep{Plexos}. Another example is the COMPETES-TNO model \citep{MoralesEspaña2024} that has incorporated hydrogen sector-specific constraints to a power system model.  The trend has been driven by the need to model decarbonization pathways involving several sector coupling technologies, e.g., hydrogen production by electrolysis, electric vehicles, and heat production from electricity \citep{Gea-Bermudez2021}. The main drawback of the process-specific approach is that the number of possible interactions between energy system components can become very large, making the model unwieldy to maintain and expand. Besides, there are only a limited number of ways processes can be described in linear programmes (LPs) or mixed-integer linear programmes (MIPs). 

Another strategy is to formulate methods between higher-level concepts like nodes, connections and units that act as building blocks (BB) to build a more general/generic energy system. Each BB can then choose an appropriate method for every particular process. This limits the number of formulations to the number of supported methods; however, it provides flexibility to the user on the modelling options. It also means that users can add new process types without writing new code, provided an appropriate method is available in the model. This approach has been used by models like the IRENA FlexTool \citep{FlexTool} and, to an extent, SpineOpt \citep{Ihlemann2022}.

A further perspective is that energy system models can be seen as network graphs that illustrate the connections among various energy assets in different sectors or energy carriers \citep{Markensteijn2020}. From this viewpoint, different types of \textit{nodes} typically serve as the basic building blocks that link all the elements, such as producers, converters, storage, and consumers. However, the existing literature has not delved into the potential use of the natural graph structure of energy systems to link the assets directly in order to improve the computational efficiency of models.

\FloatBarrier
\subsection{Contribution}
\label{sec:contribution}
The main contributions of this paper are twofold:
\begin{enumerate}
    \item In this paper, we debunk the misconception that an LP formulation cannot be further improved (sped up) without sacrificing its fidelity. We show four different approaches using different building blocks and compare the computational performance of three of them. Although all the different formulations lead to the same optimal results, they greatly differ in their computational performance, thus demonstrating that the quality of an LP model can be improved while retaining its fidelity.

     \item  We propose the \textit{Energy Asset} as a single building block, allowing a more direct connection between components by leveraging the graph structure of energy systems. The proposed strategy  replaces the traditional building block, nodes, units and connections with only one BB: \textit{Energy Assets} (vertices) and use energy flows (arcs) to connect them.  Thus, it inherently avoids unnecessary elements in between with their associated extra constraints and variables,  while keeping the full model flexible and speeding up solving times without sacrificing its fidelity.
\end{enumerate}

\section{Quality of LP models}
\label{sec:LPquality}
A common misconception is that an LP is the simplest representation of a problem for a desired fidelity, and that the model cannot be simplified without compromising the value of its solution \citep{KLOTZ20131}.
That is, the only way to speed up an LP without losing accuracy is through improvements in computing power (hardware) and LP solvers (software). This belief can lead modellers to inadvertently create large overhead model sizes, assuming the model is very efficient since LP is the most simplified you can go. When the model size is still potentially problematic, the common belief is that the only other option to speed up large-scale LPs is to solve an (over)simplified, smaller model, which sacrifices its fidelity. 

However, all models are not the same even if they model exactly the same problem (i.e., same model fidelity). One model can be faster than another under the same hardware and software. They differ in their theoretical model quality, and their quality can be improved so they can solve faster. Improving the quality of a model means that the model is reformulated so that it retains its fidelity while solving faster (using the same software and hardware). Crafting high-quality formulations allows us to increase the model fidelity without increasing solving times or even create higher quality models that solve faster, thus pushing the Pareto front model fidelity vs computational burden.

Although discussing model quality in LP Models is not common, model quality is a well-known concept in mixed-integer programming (MIP). The quality of an MIP model is defined by its tightness, that is, how near is its relaxed LP feasible region to that of the integer one. The tightest possible model (convex hull) can solve an MIP as an LP, greatly lowering the computational burden. However, trying to tighten an MIP formulation often implies increasing its size; hence, there is a trade-off between the tightness and compactness of an MIP.

What then defines the model quality of an LP model? The quality of an LP model is defined by its size. That is, a more compact model, i.e., fewer constraints/variables and non-zeros, has higher quality than a less compact one with the same fidelity, and hence, it is expected to solve faster.
Here, we differentiate the model quality with numerical-related issues, which can also slow down solving times. That is, model quality is independent of the data used. Of course, the modeller should be careful when populating the model with data to avoid LP numerical issues, such as degeneracy, numerical stability and ill-conditioning \citep{KLOTZ20131}.

How can the quality of an LP model be improved? How can it be lowered in size without sacrificing its fidelity? To the best of our knowledge, there are three ways to lower one of the dimensions of a problem. First, the trivial option is to remove equalities. Each equality means that a variable in that equality can be removed by replacing its equivalent in the other constraints. Although this procedure lowers the number of variables, it increases the number of non-zeros. Second, the Fourier–Motzkin procedure eliminates a set of variables, creating another model in which both models have the same solutions over the remaining variables. Fourier–Motzkin procedure could also eliminate constraints if applied to the dual formulation. However, this procedure comes at the expense of producing an often exponential number of constraints and non-zeros, thus creating worse bottlenecks and slowing down LP-solving times. The third way of LP reformulation involves splitting dense columns into sparser ones \citep{Irvin_1991}. Although this procedure increases the number of constraints and variables, it could speed up solving times when the number of non-zeros is the bottleneck. 

Lowering the size of an LP model is not a trivial task since lowering one dimension comes at the expense of sacrificing another dimension, which can become a new bottleneck, potentially damaging the quality of the model instead of improving it. In this paper, we emphasise in a change of paradigm for flexible energy system modelling; we use a single building block, the energy asset, thus fully exploiting the network nature of the system. This reformulation simultaneously lowers all three dimensions of the problem: including the number of variables, constraints,, and non-zeros. This higher-quality LP reformulation naturally lowers both model creation and solving times, obtaining more significant advantages as the problem size increases.

\section{Energy systems modelling and building blocks for flexible modelling} 
\label{sec:method}
Energy systems can be represented as directed graphs, where the \textit{vertices} often denote balancing nodes, and the \textit{edges} represent the energy flows between these components. In this paper, we denote the graph's vertices as the energy system's Building Blocks (BB). Depending on their characteristics in the energy system, these building blocks represent \textit{Energy Assets} that can produce, consume, convert, store, or balance energy.

Current state-of-the-art energy systems models use at least two or three building blocks, such as nodes, connections and units. Energy system models often use \textit{nodes} as an additional building block to connect energy assets and establish energy balances \citep{PyPSA}. Some models have a storage option in their node balance \citep{FlexTool}, while others maintain a separate building block for it \citep{OSeMOSYS}. Moreover, some models include an extra building block to describe the \textit{connections} between nodes \citep{Ihlemann2022}. Although having multiple building blocks to model an energy system seems appealing due to the flexibility to represent various configurations, it often comes with a computational cost that becomes more significant for large-scale problems.

In this paper, we highlight the advantages of using only one building block: energy assets, represented as the vertices of a graph, with the flows between the assets as the arcs. Adopting this approach makes the node concept irrelevant since we directly connect assets. This approach simplifies the graph structure of the energy system, allowing for flexible modelling options that reduce the number of variables and constraints required to model the same problem. This reduction has computational benefits in both creating the optimisation model and the time to solve it, as demonstrated in the experiments discussed in Section \ref{sec:results}.

\subsection{Generic formulation} 
\label{sec:generic-formulation}

While each approach may have more specific constraints and differing notations, we can generalize them all by using  the proposed \textit{energy asset} as the only building block. The following formulation makes comparing the formulations in the following sections more straightforward.

\begin{equation}
\label{eq:objective_function}
\begin{aligned}
\text{{minimize}} \quad \sum_{a \in \mathcal{A}^{\text{i}}} C^{I}_{a} \cdot i_{a} + \sum_{(a^{\text{from}},a^{\text{to}}) \in \mathcal{F}} \sum_{t \in \mathcal{T}} C^{O}_{(a^{\text{from}},a^{\text{to}})} \cdot f_{(a^{\text{from}},a^{\text{to}}),t}
\end{aligned}
\end{equation} 

$s.t.$

\begin{equation}
\label{eq:consumer_balance}
\begin{aligned}
\sum_{(a^{\text{from}},a) \in \mathcal{F}} f_{(a^{\text{from}},a),t} - \sum_{(a,a^{\text{to}}) \in \mathcal{F}} f_{(a,a^{\text{to}}),t} = D_{a,t} \quad \forall a \in \mathcal{A}^{\text{c}}, \forall t \in \mathcal{T}
\end{aligned}
\end{equation}

\begin{equation}
\label{eq:storage_balance}
\begin{aligned}
s_{a,t} = & s_{a,t-1 | t>1} + S^{0}_{a | t=1} \\
& + \eta^{\text{in}}_{\text{a}} \cdot \sum_{(a^{\text{from}},a) \in \mathcal{F}} f_{(a^{\text{from}},a),t} \\
& - \frac{1}{\eta^{\text{out}}_{\text{a}}} \cdot \sum_{(a,a^{\text{to}}) \in \mathcal{F}} f_{(a,a^{\text{to}}),t}  \quad \forall a \in \mathcal{A}^{\text{s}}, \forall t \in \mathcal{T}
\end{aligned}
\end{equation}

\begin{equation}
\label{eq:conversion_balance}
\begin{aligned}
\eta^{\text{in}}_{\text{a}} \cdot \sum_{(a^{\text{from}},a) \in \mathcal{F}} f_{(a^{\text{from}},a),t} = \sum_{(a,a^{\text{to}}) \in \mathcal{F}} f_{(a,a^{\text{to}}),t}  \quad \forall a \in \mathcal{A}^{\text{cv}}, \forall t \in \mathcal{T}
\end{aligned}
\end{equation}

\begin{equation}
\label{eq:capacity_limits}
\begin{aligned}
\sum_{(a,a^{\text{to}}) \in \mathcal{F}} f_{(a,a^{\text{to}}),t} \leq \overline{P}_a \cdot (U^0_{a} + i_{a}) \quad
\forall a \in \mathcal{A}, \forall t \in \mathcal{T}
\end{aligned}
\end{equation}

\begin{equation}
\label{eq:charging_limits}
\begin{aligned}
\sum_{(a^{\text{from}},a) \in \mathcal{F}} f_{(a^{\text{from}},a),t} \leq \overline{P}_a \cdot (U^0_{a} + i_{a}) \quad
\forall a \in \mathcal{A}^{\text{s}}, \forall t \in \mathcal{T}
\end{aligned}
\end{equation}

\begin{equation}
\label{eq:storage_capacity}
\begin{aligned}
0 \leq s_{a,t} \leq \overline{S}_{a} \quad
\forall a \in \mathcal{A}^{\text{s}}, \forall t \in \mathcal{T}
\end{aligned}
\end{equation}

\begin{equation}
\label{eq:flow_bounds}
\begin{aligned}
-\overline{F}_{(a^{\text{to}},a^{\text{from}})} \leq f_{(a^{\text{from}},a^{\text{to}}),t} \leq \overline{F}_{(a^{\text{from}},a^{\text{to}})} \quad
\forall (a^{\text{from}},a^{\text{to}}) \in \mathcal{F}, \forall t \in \mathcal{T}
\end{aligned}
\end{equation}

\begin{equation}
\label{eq:invest_limit}
\begin{aligned}
i_{a} \leq \overline{I}_{a} \quad
\forall  a \in \mathcal{A}^{\text{inv}}
\end{aligned}
\end{equation}

\begin{equation}
\label{eq:no_invest}
\begin{aligned}
i_{a} = 0 \quad
\forall  a \notin \mathcal{A}^{\text{inv}}
\end{aligned}
\end{equation}

Where:
\begin{itemize}
    \item[] $\mathcal{A}$: set of all assets with elements $a \in \mathcal{A}$
    
    \item[] $\mathcal{A}^{\text{i}}  \subseteq \mathcal{A}$: subset of investable assets 
    
    \item[] $\mathcal{A}^{\text{c}}  \subseteq \mathcal{A}$: subset of consumer assets 
    \item[] $\mathcal{A}^{\text{s}}  \subseteq \mathcal{A}$: subset of storage assets     
    \item[] $\mathcal{A}^{\text{cv}} \subseteq \mathcal{A}$: subset of conversion assets     
    \item[] $\mathcal{F}$: set of flows between two assets with elements $(a^{\text{from}},a^{\text{to}}) \in \mathcal{F}$
    \item[] $\mathcal{T}$: set of all timesteps with elements $t \in \mathcal{T}$

    \item[] $i_{a}$: investment unit variable of asset $a$

    \item[] $f_{(a^{\text{from}},a^{\text{to}}),t}$: flow variable from asset $a^{\text{from}}$ to asset $a^{\text{to}}$ in timestep $t$
    \item[] $s_{a,t}$: storage level variable of asset $a$ in timestep $t$    
    \item[] $C^{I}_{a}$: investment cost parameter of asset $a$
    \item[] $C^{O}_{(a^{\text{from}},a^{\text{to}})}$: operation cost parameter of flow between two assets
    \item[] $D_{a,t}$: demand parameter of the consumer asset $a$ in timestep $t$
    \item[] $S^0_{a}$: initial storage level parameter of the storage asset $a$

    \item[] $U^0_{a}$: initial number of units parameter of the asset $a$
    \item[] $\overline{P}_{a}$: maximum capacity parameter of the asset $a$    
    \item[] $\overline{F}_{(a^{\text{from}},a^{\text{to}})}$, $\overline{F}_{(a^{\text{to}},a^{\text{from}})}$: maximum flow between two assets in both directions.    
    \item[] $\overline{S}_{a}$: maximum storage capacity parameter of the storage asset $a$
    \item[] $\overline{I}_{a}$: maximum investment limit parameter of the asset $a$  

    \item[] $\eta^{\text{in}}_{\text{a}}$: efficiency parameter of flows going into the asset $a$
    \item[] $\eta^{\text{out}}_{\text{a}}$: efficiency parameter of flows going out the asset $a$
\end{itemize}

The objective function (\ref{eq:objective_function}) minimises the 
 investment cost and the operation 
variable cost of the flow between two assets. Equations (\ref{eq:consumer_balance}), (\ref{eq:storage_balance}) and (\ref{eq:conversion_balance}) define the balance constraint for consumer, storage, and conversion assets. The capacity limit of the flows (\ref{eq:capacity_limits}) is defined by the maximum capacity of the asset, while the maximum charging limit for storage assets is defined by  (\ref{eq:charging_limits}). Equation (\ref{eq:storage_capacity}) represents the maximum storage capacity, and  (\ref{eq:flow_bounds}) defines the bounds for the flow variable.  Equation (\ref{eq:invest_limit}) and (\ref{eq:no_invest}) limit the availability of the investment variables for the assets. 
Annex A shows how to model a DC power flow and unit commitment in the proposed approach. In addition, the TulipaEnergyModel \cite{TulipaEnergyModel2023} documentation\footnote{\url{https://tulipaenergy.github.io/TulipaEnergyModel.jl/dev/40-formulation/}} shows a comprehensive overview of the proposed approach.

\subsection{Illustrative example} 
\label{sec:illustrative-example}

\begin{figure}
    \centering
    \includegraphics[width=\textwidth]{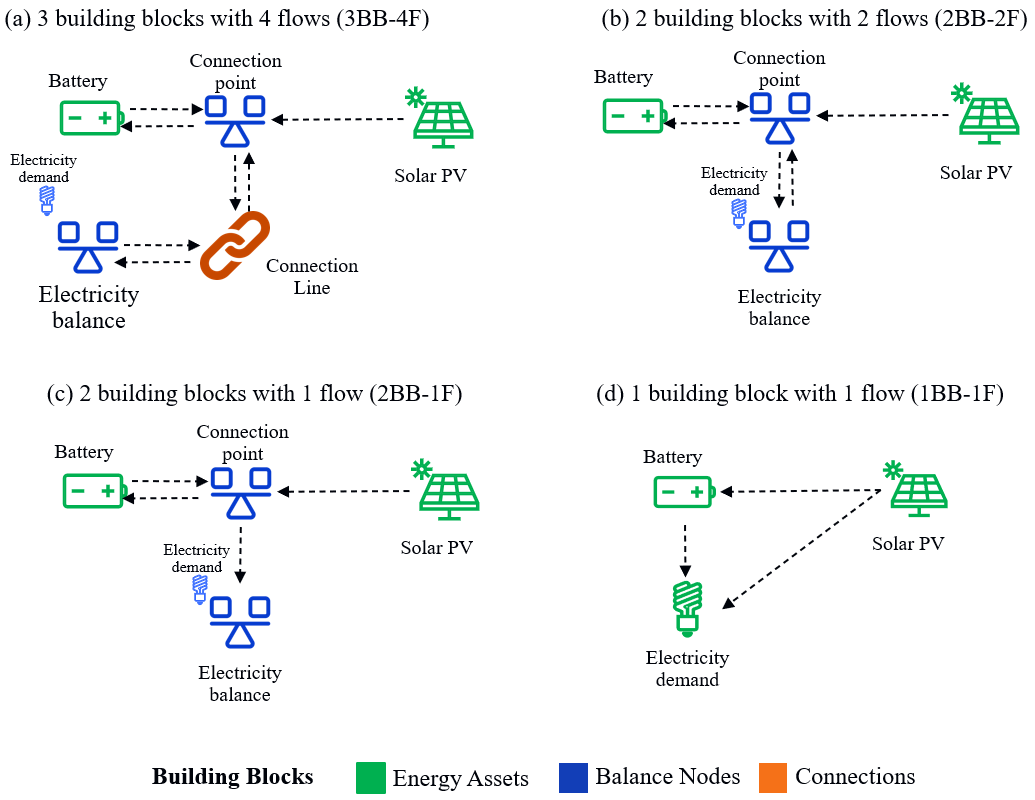}
    \caption{Using different BB to model a hybrid solar PV with storage. In general, every arrow represents one variable, and every energy asset one constraint}
    \label{fig:example-different-modelling}
\end{figure}

Having flexible modelling options is becoming more relevant in energy system analysis with the appearance of new configurations, such as the hybrid configurations of storage and renewable assets. Figure \ref{fig:example-different-modelling} shows an example of such a hybrid configuration, where storage can only be charged from solar PV. Still, the storage and solar PV assets can both deliver energy to the grid. We can take this use case as an example to understand the difference between various approaches,  which only considers operation variables for the sake of simplicity . In this case, the storage asset can only be charged from the renewable asset and not from the grid.  Figure \ref{fig:example-different-modelling} gives an overview of each modelling approach for this example. In the following sections, we describe each approach in detail. It is worth noting that different modelling approaches represent the same situation, i.e., they are all modelling exactly the same problem, and this example is not unique. For instance, green hydrogen production is another example where a conversion asset, such as an electrolyser, can only produce hydrogen from the renewable asset; however, the renewable asset can still send energy to the network.

\FloatBarrier
\subsubsection{Three building blocks with four flow variables (3BB-4F)} 
\label{sec:nodal-4-flows}
This approach uses three BB:  energy \textit{assets},  balance \textit{nodes}, and  \textit{connections} to link the nodes. In addition, the connection uses four positive variables to link the nodes, two for each node it connects to. This approach is advantageous as it can represent situations where the incoming flow might differ from the outgoing flow, which is needed to model transmission or distribution network losses or gas flows with linepack. However, modelling network constraints without losses or differences between the incoming and outgoing flow creates extra variables and constraints.

Figure \ref{fig:example-different-modelling}(a) shows the battery (bt) with a renewable (pv) example for the 2BB-4F. The connection point (cp) is an auxiliary node connecting the two assets. The electricity demand (ed) is balanced in an extra node. Moreover, the maximum capacity of the flow coming from the connection line (cl) to the connection point must be zero to avoid charging the battery from the grid (i.e., electricity demand balance node). Therefore, we need 8 variables (7 flow variables $f$ + 1 storage level $s$) and 11 constraints (8 capacity limits + 2 node balances + 1 storage balance). The index $t$ is not included in the following equations for simplicity. 

\begin{itemize}
    \item Connection point node balance constraint:
\begin{align*}
  f_{\text{(bt,cp)}} + f_{\text{(pv,cp)}} + f_{\text{(cl,cp)}} = f_{\text{(cp,bt)}} + f_{\text{(cp,cl)}}
\end{align*}
    \item Electricity demand node balance constraint:
\begin{align*}
  f_{\text{(cl,ed)}} = D_{\text{ed}} + f_{\text{(ed,cl)}}
\end{align*}
    \item Battery storage balance constraint:
\begin{align*}
  s_{\text{bt}} = S^{0}_{\text{bt}} + \eta^{\text{in}}_{\text{bt}} \cdot f_{\text{(cp,bt)}} - \frac{1}{\eta^{\text{out}}_{\text{bt}}} \cdot f_{\text{(bt,cp)}}
\end{align*}
    \item Capacity limit constraints:
\begin{align*}
&f_{\text{(bt,cp)}} \leq P_{\text{bt}}, \quad \quad  
 f_{\text{(pv,cp)}} \leq P_{\text{pv}}, \quad \quad
 f_{\text{(cl,cp)}} \leq F_{\text{(cl,cp)}}, \\
&f_{\text{(cp,bt)}} \leq P_{\text{bt}}, \quad \quad 
 f_{\text{(cp,cl)}} \leq F_{\text{(cp,cl)}}, \quad \quad
 f_{\text{(cl,ed)}} \leq F_{\text{(cl,ed)}}, \\ 
&f_{\text{(ed,cl)}} \leq F_{\text{(ed,cl)}}, \quad \quad
 s_{\text{bt}} \leq S_{\text{bt}} \quad \quad
\end{align*}
    \item Non-negative variables:
\begin{align*}
f_{\text{(bt,cp)}},
f_{\text{(pv,cp)}},
f_{\text{(cl,cp)}},
f_{\text{(cp,bt)}},
f_{\text{(cp,cl)}},
f_{\text{(cl,ed)}},
f_{\text{(ed,cl)}},
s_{\text{bt}} \geq 0
\end{align*}
\end{itemize}

Where:

\begin{itemize}
    \item[] $f_{\text{(bt,cp)}}$: flow from the battery to the connection point (i.e., battery discharge) limited by battery capacity $P_{\text{bt}}$
    \item[] $f_{\text{(cp,bt)}}$: flow from the connection point to the battery (i.e., battery charge) limited by battery capacity $P_{\text{bt}}$\footnote{Charging and discharging capacities may vary based on the type of storage asset; for simplicity, we consider them equal.}   
    \item[] $f_{\text{(pv,cp)}}$: flow from the solar pv to the connection point limited by solar PV availability $P_{\text{pv}}$       
    \item[] $f_{\text{(cl,cp)}}$: flow from the connection line to the connection point with a maximum capacity $F_{\text{(cl,cp)}}$           
    \item[] $f_{\text{(cp,cl)}}$: flow from the connection point to the connection line with a maximum capacity $F_{\text{(cp,cl)}}$               
    \item[] $f_{\text{(cl,ed)}}$: flow from the connection line to the electricity demand balance with a maximum capacity $F_{\text{(cl,ed)}}$                   
    \item[] $f_{\text{(ed,cl)}}$: flow from the electricity demand balance to the connection line with a maximum capacity $F_{\text{(ed,cl)}}$
    \item[] $s_{\text{bt}}$: storage level of the battery with a maximum capacity $S_{\text{bt}}$   
    \item[] $D_{\text{ed}}$: electricity demand input data   
    \item[] $\eta^{\text{in}}_{\text{bt}}$ and $\eta^{\text{out}}_{\text{bt}}$: charging and discharging efficiencies of the battery   
    \item[] $S^{0}_{\text{bt}}$: initial storage level of the battery       
\end{itemize}

\FloatBarrier
\subsubsection{Two building blocks with two flow variables (2BB-2F)} 
\label{sec:nodal-2-flows}
This approach uses the energy \textit{assets} and the balance \textit{nodes} as building blocks. In addition, it uses two positive flow variables to connect the nodes instead of the four variables used in the previous approach. The two flow variables are needed to model transmission or distribution network losses. However, gas flows with linepack cannot be modelled correctly with only two flow variables. Depending on the input data, the incoming and outgoing flow can have values simultaneously, especially when losses are not considered. In such cases, modellers typically use the net value as the transfer between the nodes. Alternatively, they may include a binary variable to prevent simultaneous incoming and outgoing flows, although this would convert the problem into a Mixed-Integer Programming problem.

Figure \ref{fig:example-different-modelling}(b) shows the battery with a renewable example for the 2BB-2F. The connection point is an auxiliary node connecting the two assets. Moreover, the maximum capacity of the flow coming from the electricity balance node to the connection point must be zero to avoid charging the battery from the grid. Therefore, we need 6 variables (5 flow variables + 1 storage level) and 9 constraints (6 capacity limits + 2 node balances + 1 storage balance). 

\begin{itemize}
    \item Connection point node balance constraint:
\begin{align*}
  f_{\text{(bt,cp)}} + f_{\text{(pv,cp)}} + f_{\text{(ed,cp)}} = f_{\text{(cp,bt)}} + f_{\text{(cp,ed)}}
\end{align*}
    \item Electricity demand node balance constraint:
\begin{align*}
  f_{\text{(cp,ed)}} = D_{\text{ed}} + f_{\text{(ed,cp)}}
\end{align*}
    \item Battery storage balance constraint:
\begin{align*}
  s_{\text{bt}} = S^{0}_{\text{bt}} + \eta^{\text{in}}_{\text{bt}} \cdot f_{\text{(cp,bt)}} - \frac{1}{\eta^{\text{out}}_{\text{bt}}} \cdot f_{\text{(bt,cp)}}
\end{align*}
    \item Capacity limit constraints:
\begin{align*}
&f_{\text{(bt,cp)}} \leq P_{\text{bt}}, \quad \quad  
 f_{\text{(pv,cp)}} \leq P_{\text{pv}}, \quad \quad
 f_{\text{(ed,cp)}} \leq F_{\text{(ed,cp)}}, \\
&f_{\text{(cp,bt)}} \leq P_{\text{bt}}, \quad \quad 
 f_{\text{(cp,ed)}} \leq F_{\text{(cp,ed)}}, \quad \quad
 s_{\text{bt}} \leq S_{\text{bt}} \quad \quad
\end{align*}
    \item Non-negative variables:
\begin{align*}
f_{\text{(bt,cp)}},
f_{\text{(pv,cp)}},
f_{\text{(ed,cp)}},
f_{\text{(cp,bt)}},
f_{\text{(cp,ed)}},
s_{\text{bt}} \geq 0
\end{align*}
\end{itemize}

Where the new variables are:

\begin{itemize}
    \item[] $f_{\text{(cp,ed)}}$: flow from the connection point to the electricity demand balance with a maximum capacity $F_{\text{(cp,ed)}}$                   
    \item[] $f_{\text{(ed,cp)}}$: flow from the electricity demand balance to the connection point with a maximum capacity $F_{\text{(ed,cp)}}$      
\end{itemize}

\FloatBarrier
\subsubsection{Two building blocks with one flow variable (2BB-1F)} 
\label{sec:nodal-1-flow}
This approach uses the energy \textit{assets} and the balance \textit{nodes} as building blocks. In addition, it uses a single free variable to represent the flow between nodes, which can take positive and negative values instead of two positive variables as in the previous approach. The main advantage of this method is that it eliminates the possibility of bidirectional flow between two nodes, as there is only one variable. However, it is not possible to model transmission or distribution network losses using a single free variable for the flow. Lastly, the free variable must have bounds in both directions to represent the capacity limits between the two nodes.

Figure \ref{fig:example-different-modelling}(c) shows the battery with a renewable example for this approach. The connection point is an auxiliary node connecting the two assets. Moreover, the maximum capacity of the flow in the direction from the electricity balance node to the connection point must be zero to avoid charging the battery from the grid. Therefore, we need 5 variables (4 flow variables + 1 storage level) and 9 constraints (6 capacity limits + 2 node balances + 1 storage balance). 

\begin{itemize}
    \item Connection point node balance constraint:
\begin{align*}
  f_{\text{(bt,cp)}} + f_{\text{(pv,cp)}} = f_{\text{(cp,bt)}} + f_{\text{(cp,ed)}}
\end{align*}
    \item Electricity demand node balance constraint:
\begin{align*}
  f_{\text{(cp,ed)}} = D_{\text{ed}}
\end{align*}
    \item Battery storage balance constraint:
\begin{align*}
  s_{\text{bt}} = S^{0}_{\text{bt}} + \eta^{\text{in}}_{\text{bt}} \cdot f_{\text{(cp,bt)}} - \frac{1}{\eta^{\text{out}}_{\text{bt}}} \cdot f_{\text{(bt,cp)}}
\end{align*}
    \item Capacity limit constraints:
\begin{align*}
&f_{\text{(bt,cp)}} \leq P_{\text{bt}}, \quad \quad  
 f_{\text{(pv,cp)}} \leq P_{\text{pv}}, \quad \quad
 f_{\text{(cp,bt)}} \leq P_{\text{(cp,bt)}}, \\
& -F_{\text{(ed,cp)}} \leq f_{\text{(cp,ed)}} \leq F_{\text{(cp,ed)}}, \quad \quad
 s_{\text{bt}} \leq S_{\text{bt}} \quad \quad
\end{align*}
    \item Non-negative variables:
\begin{align*}
f_{\text{(bt,cp)}},
f_{\text{(pv,cp)}},
f_{\text{(cp,bt)}},
s_{\text{bt}} \geq 0
\end{align*}
\end{itemize}

\FloatBarrier
\subsubsection{One building block with one flow variable (1BB-1F)} 
\label{sec:asset-to-asset}
This approach only uses the energy \textit{assets} as building blocks. Using the graph-theory principles establishes the connection between energy assets as vertices and energy flows as edges. Connecting assets directly to each other (without any intervening nodes) can significantly reduce the number of variables and constraints required to represent the system. When transmission or distribution network losses are significant, two conversion assets can represent them—one for each direction with input-output ratios. A similar approach can be adopted for gas pressure flows. Therefore, this method enables us to easily model simple and complex situations, reducing the model size while retaining the same accuracy.

Figure \ref{fig:example-different-modelling}(d) shows the battery with a renewable example for the 1BB-1F approach. Since the assets can connect among them, the battery can directly charge from renewable, and there is no need for an extra constraint to avoid charging from the grid. Therefore, we need 4 variables (3 flow variables + 1 storage level) and 6 constraints (4 capacity limits + 1 demand balance + 1 storage balance). 

\begin{itemize}
    \item Electricity demand node balance constraint:
\begin{align*}
  f_{\text{(bt,ed)}} + f_{\text{(pv,ed)}} = D_{\text{ed}} 
\end{align*}
    \item Battery storage balance constraint:
\begin{align*}
  s_{\text{bt}} = S^{0}_{\text{bt}} + \eta^{\text{in}}_{\text{bt}} \cdot f_{\text{(pv,bt)}} - \frac{1}{\eta^{\text{out}}_{\text{bt}}} \cdot f_{\text{(bt,ed)}}
\end{align*}
    \item Capacity limit constraints:
\begin{align*}
&f_{\text{(bt,ed)}} \leq P_{\text{bt}}, \quad \quad  
 f_{\text{(pv,bt)}} + f_{\text{(pv,ed)}} \leq P_{\text{pv}}, \quad \quad  \\
&f_{\text{(pv,bt)}} \leq P_{\text{bt}}, \quad \quad   
 s_{\text{bt}} \leq S_{\text{bt}} \quad \quad
\end{align*}
    \item Non-negative variables:
\begin{align*}
f_{\text{(bt,ed)}},
f_{\text{(pv,bt)}},
f_{\text{(pv,ed)}},
s_{\text{bt}} \geq 0
\end{align*}
\end{itemize}

Where the new variables are:

\begin{itemize}
    \item[] $f_{\text{(bt,ed)}}$: flow from the battery to the electricity demand balance limited by battery capacity $P_{\text{bt}}$                   
    \item[] $f_{\text{(pv,bt)}}$: flow from the solar pv to the battery limited by battery capacity $P_{\text{bt}}$ and the solar PV availability $P_{\text{pv}}$
    \item[] $f_{\text{(pv,ed)}}$: flow from the solar pv to the electricity demand balance limited by the solar availability $P_{\text{pv}}$
\end{itemize}

\FloatBarrier
\subsubsection{Summary} 
Table \ref{tab:comparison-tiny-case} summarises the number of variables and constraints per time step for each modelling approach in the example shown in Figure \ref{fig:example-different-modelling}. It also shows the reduction in variables and constraints, with the 3BB-4F approach as a reference.

Reducing the number of variables and constraints significantly benefits the time taken to build and solve an optimisation problem, as we show in Section \ref{sec:results}. Although the solvers' presolve can eliminate unnecessary variables and constraints, we show that modellers can further improve this process by using formulations with fewer variables and constraints while representing the same energy system. Thus speeding up solving times.

\begin{table}[ht]
\centering
\caption{Number of variables and constraints in each modelling approach per time step}
\label{tab:comparison-tiny-case}
\begin{tabular}{ccccc}
\toprule
\textbf{Modelling Approach} & \textbf{Variables} & \textbf{Constraints} & \textbf{Non-zeros} \\ \hline
3BB-4F  &    8	&  11	& 18\\	
2BB-2F  &    6	($\downarrow$ 25\%)&  9	($\downarrow$ 18\%) &  16	($\downarrow$ 11\%)\\	
2BB-1F  &    5	($\downarrow$ 38\%)&  9	($\downarrow$ 18\%) &  13	($\downarrow$ 28\%)\\		
1BB-1F	&    4	($\downarrow$ 50\%)&  6	($\downarrow$ 45\%) &  9	($\downarrow$ 50\%)\\ \bottomrule
\end{tabular}
\end{table}

\FloatBarrier
\section{Case Study} 
\label{sec:calculation}

The following sections describe the energy system optimisation model and a case study that compares different approaches. For the optimisation model, we have selected TulipaEnergyModel.jl \citep{TulipaEnergyModel2023} as it can model all the approaches based on the input data. The case study highlights situations where the 1BB-1F approach can leverage its flexibility to connect energy assets; nevertheless, Section \ref{sec:discussion} discusses when this is possible and provides some insights for energy system modellers.
\FloatBarrier
\subsection{Energy system optimisation model} 
\label{sec:model}
TulipaEnergyModel.jl is an optimisation model using 1BB-1F and determines the optimal investment and operation decisions for different types of assets (e.g., producers, consumers, conversion, storage, and transport). It is developed in Julia \citep{Bezanson_Julia_A_fresh_2017} and depends mainly on the JuMP.jl \citep{Lubin2023} and Graphs.jl \citep{Graphs2021} packages. The complete description of the model, its core concepts, mathematical formulation, and tutorials are available in the GitHub documentation of the model \footnotemark{}. 

\footnotetext{https://tulipaenergy.github.io/TulipaEnergyModel.jl/stable/}

\FloatBarrier
\subsection{Case study description} 
\label{sec:case-study}
In our case study, we conducted experiments on three approaches: 2BB-2F, 2BB-1F, and 1BB-1F. We omitted the 3BB-4F approach because the results in Section \ref{sec:results} showed it would have performed worse than the other approaches for the case study. Our focus is on illustrating the performance differences between the three selected approaches.

\begin{figure}
    \centering
    \includegraphics[width=\textwidth]{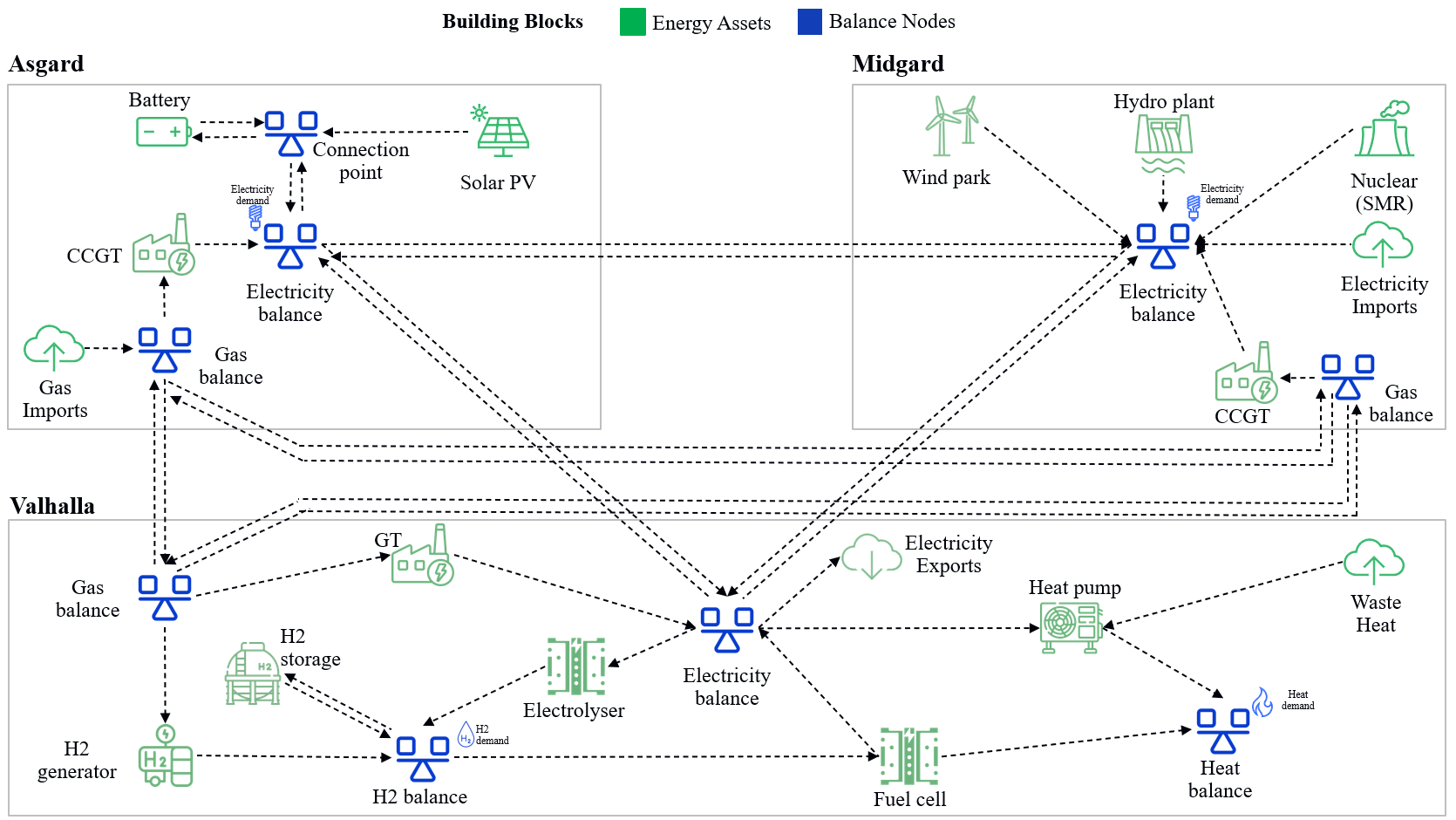}
    \caption{Case study using 2BB-2F}
    \label{fig:case-study-nodal-2-flows}
\end{figure}

Figure \ref{fig:case-study-nodal-2-flows} shows an illustrative integrated energy system with three interconnected areas: Asgard, Midgard, and Valhalla. The diagram uses the 2BB-2F approach. In addition, it includes the flow and balance of electricity, heat, and gas within a mock-up energy grid to explore different possibilities for flows among energy assets using nodes. Asgard includes a combined cycle gas turbine, a solar photovoltaic installation, and a battery system. Midgard features a wind park, a hydro plant, and a small modular reactor for nuclear power generation. Valhalla focuses on hydrogen as an energy vector, with a hydrogen generator, a hydrogen storage facility, and a fuel cell. Transmission lines and gas pipelines connect the system, allowing energy transfer between areas.

\begin{figure}
    \centering
    \includegraphics[width=\textwidth]{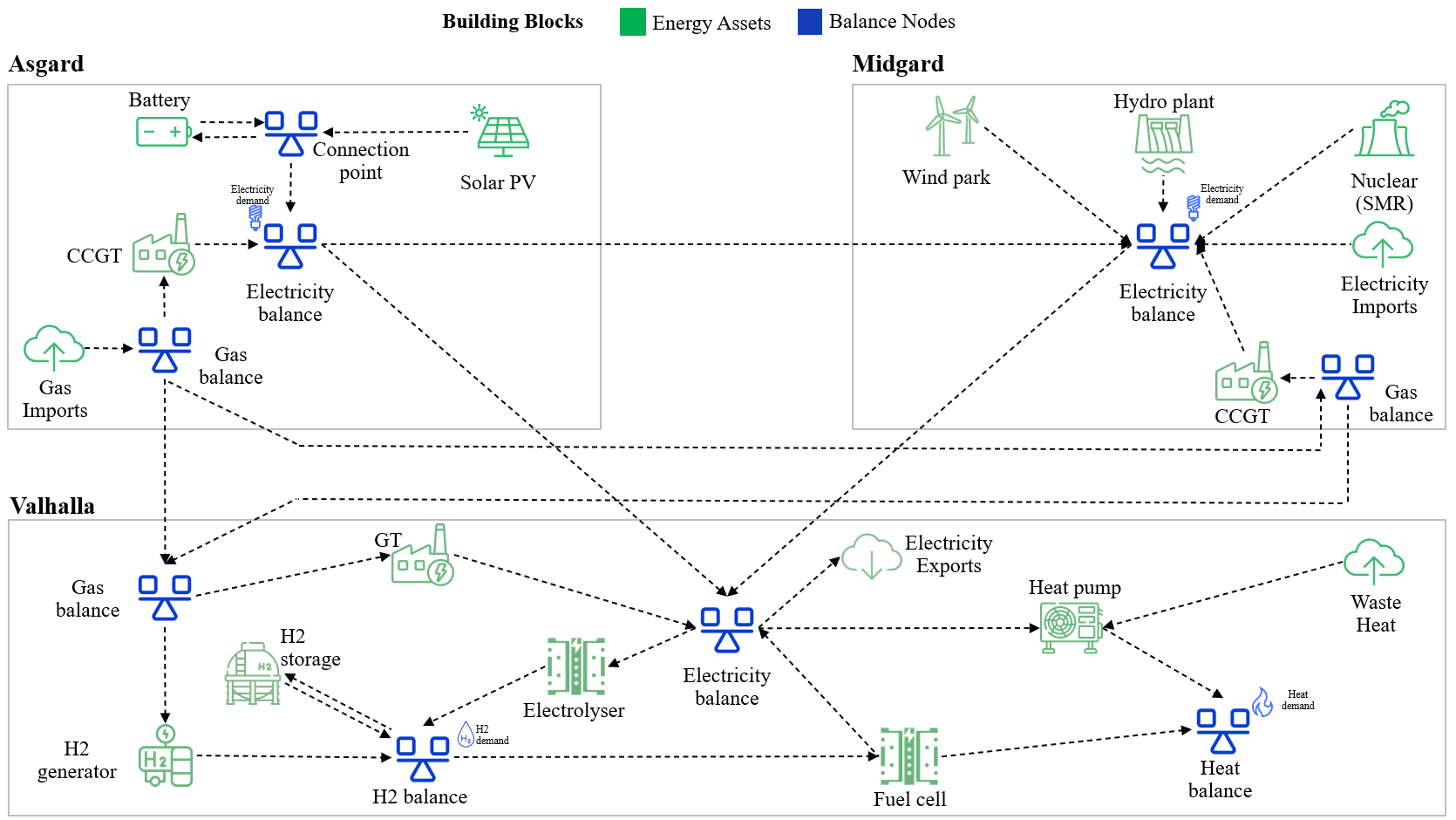}
    \caption{Case study using 2BB-1F}
    \label{fig:case-study-nodal-1-flow}
\end{figure}

\begin{figure}
    \centering
    \includegraphics[width=\textwidth]{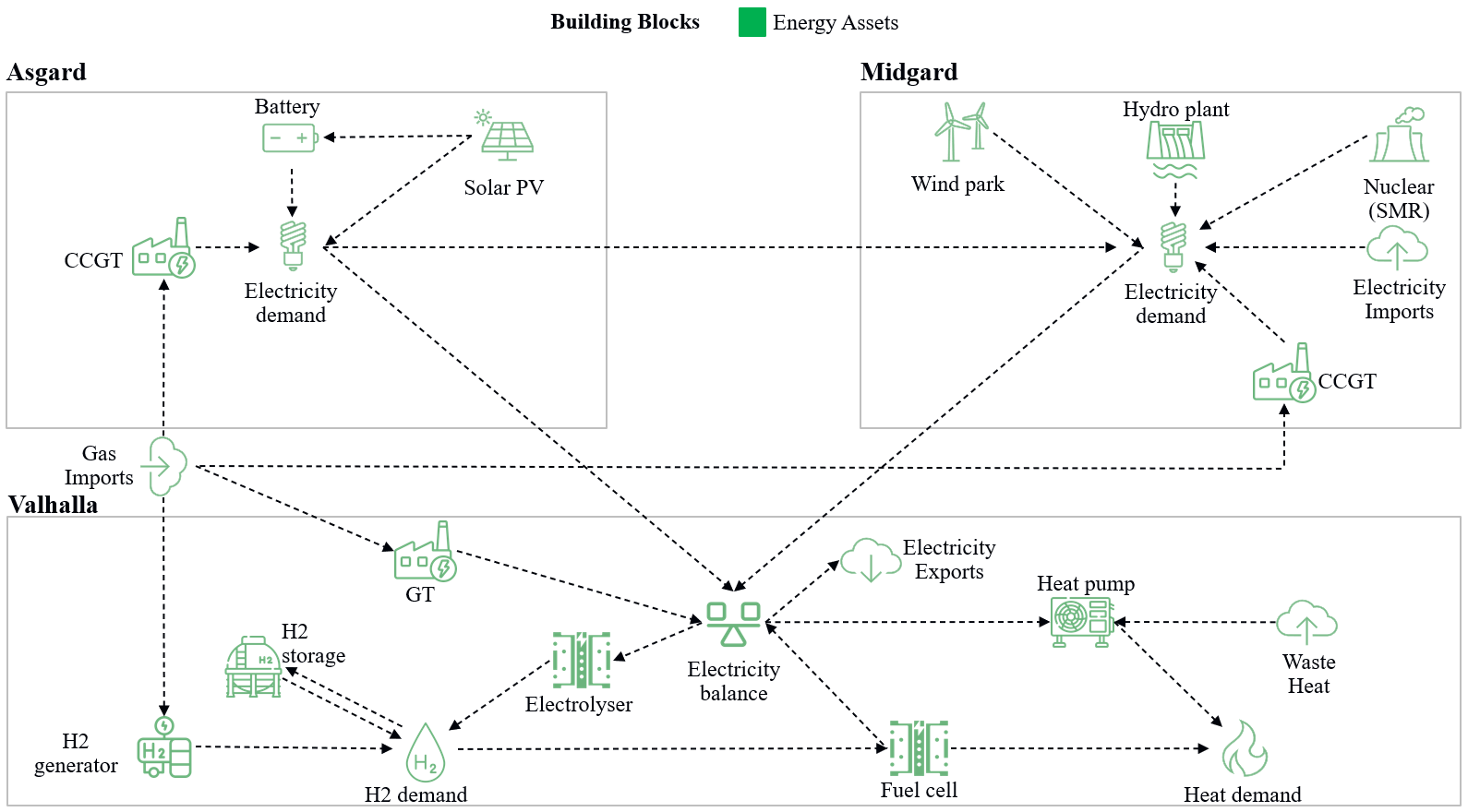}
    \caption{Case study using 1BB-1F}
    \label{fig:case-study-asset-to-asset}
\end{figure}

Figures \ref{fig:case-study-nodal-1-flow} and \ref{fig:case-study-asset-to-asset} show the equivalent energy system for the other two approaches. In general, each arrow represents a variable, and each BB element represents a constraint. Note that the number of flows (arrows) reduces compared to the 2BB-2F approach.

To assess the impact of varying problem sizes, we developed six instances of the problem, classified from small to large-scale optimisation. The smallest instance is labelled as 1, while the largest is labelled as 6. Each instance consists of hourly time steps, with differences in the time horizon covered. For example, instance 1 spans one month, whereas instance 6 covers four years. 

Section \ref{sec:results} analyses the impact of these reductions in the model for each approach and instance. Finally, Section \ref{sec:supplementary_material} has the link to the repository with the input data files for each approach.

\FloatBarrier
\section{Results} 
\label{sec:results}
The results in this section were obtained using TulipaEnergyModel.jl version 0.6.1 and Gurobi version 11.0.0 on a 12th Gen Intel(R) Core(TM) i7-1255U 1.70 GHz processor and 16.0 GB RAM. Section \ref{sec:supplementary_material} includes the link to the repository with the files to reproduce the experiments in this paper.

Table \ref{tab:comparison-sizes} shows the objective functions and problem size for each instance, where  2BB-2F is used as a reference in the table. 

We cover a wide range of model sizes, from smaller instances representing optimisation models with a few thousand variables and constraints to larger instances with millions of variables and constraints.

The reduction in the number of variables and constraints is 14\% and 18\% for the 2BB-1F approach and 26\% and 35\% for the 1BB-1F approach. It is important to highlight that the reduction is identical for all instances since the instance is an enlarged version of the same case study. Still, the connection between assets remains unchanged in all instances. It is also worth noting that the objective function is the same for every instance, regardless of the applied approach, indicating that they all represent the same energy problem. 

These results demonstrate that reducing the LP size while maintaining fidelity is possible. This finding is significant for practitioners looking to solve large-scale problems more efficiently without losing the desired level of detail.

\begin{figure}
    \centering
    \includegraphics[width=\textwidth]{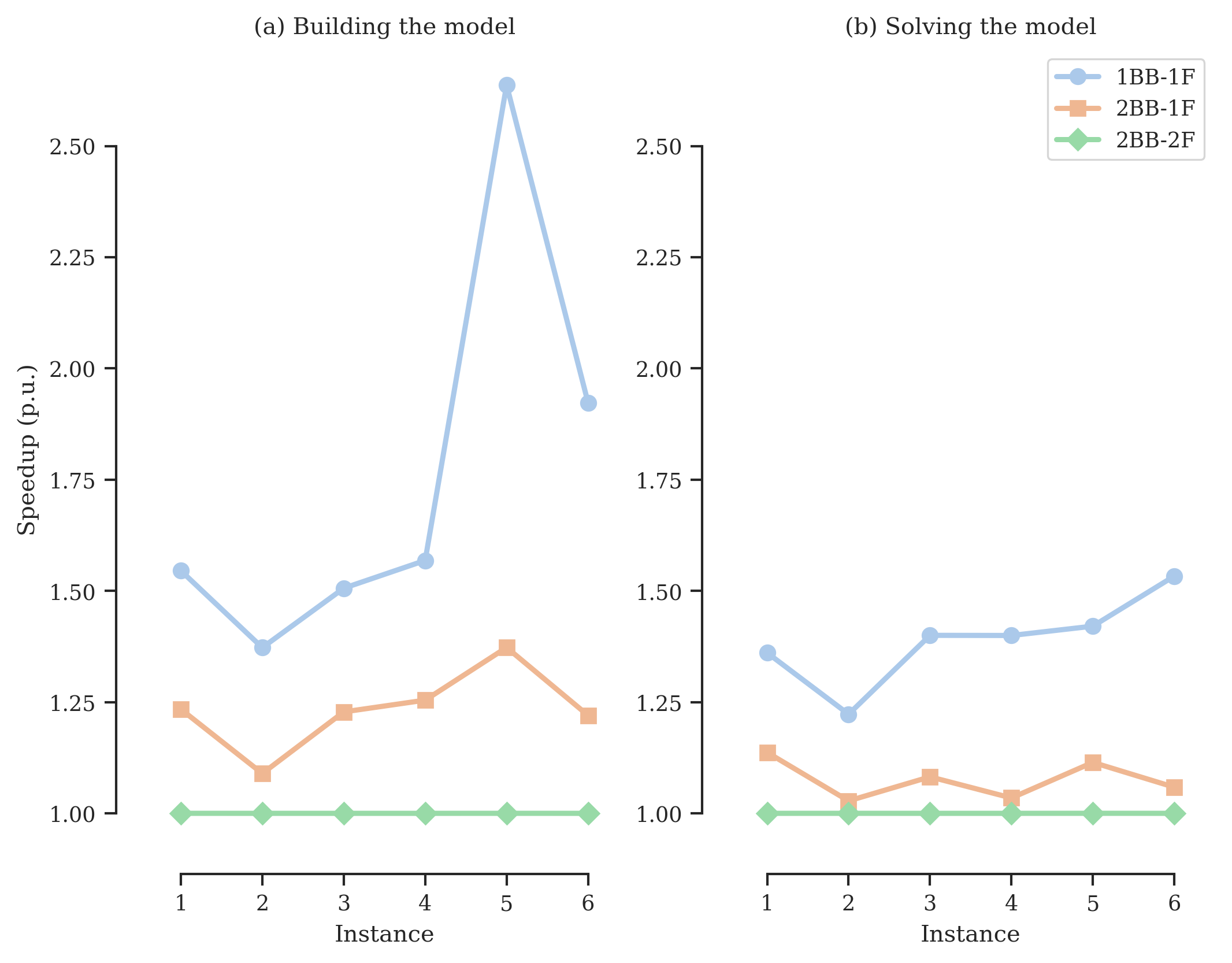}
    
    \caption{Speedups comparison}
    \label{fig:speedup-median-comparison}
\end{figure}

The graph in Figure \ref{fig:speedup-median-comparison} depicts the median speedup values for building and solving models using the three approaches and all the instances. The reference point used for comparison is the 2BB-2F approach. The values for the other approaches indicate the proportion of time the approach takes compared to the reference one. Reducing the number of variables and constraints has advantageous effects on building and solving models. Creating fewer variables and constraints has a significant advantage when building the model. Interestingly, there is a trend towards greater speedups as the instances increase. The surprising result is that having fewer variables and constraints also leads to speedups in the solving phase, with a slight trend towards increasing speedups as the instance size grows. Solvers can leverage these reductions to solve the model faster at each iteration with fewer variables and constraints. It is common practice among energy modellers to rely on the solver's presolve function to eliminate redundant variables and constraints. However, having a clean and concise formulation representing the same energy problem is always better, as it allows the solver to be faster in its process.


\begin{table}[ht]
\centering
\caption{Size of the problem in each modelling approach}
\label{tab:comparison-sizes}
\begin{tabular}{cccccc}
\toprule
\textbf{Approach} & \textbf{Instance} & \textbf{Obj. func} & \textbf{Variables} & \textbf{Constraints} & \textbf{Nonzeros}\\ \hline
         & 1	& 2.48E+08	&    28,908	&    45,696	& 96,318  \\	
         & 2 & 3.55E+08	&   173,388	&   274,176	& 578,609 \\		
2BB-2F   & 3 & 6.10E+08	 &   376,692 &   595,680 & 1,257,055 \\		
         & 4	& 1.05E+09	&   753,372	& 1,191,360	& 2,514,096 \\		
         & 5	& 1.48E+09	& 1,130,052	& 1,787,040	& 3,771,168 \\		
         & 6	& 1.89E+09	& 1,506,732	& 2,382,720	& 5,028,256 \\	\midrule	
2BB-1F   & 1-6	& 1 p.u.	&    \(\downarrow 14\%\)	&    \(\downarrow 18\%\) &    \(\downarrow 17\%\) \\ \midrule
1BB-1F   & 1-6	& 1 p.u.	&    \(\downarrow 26\%\)	&    \(\downarrow 35\%\) &    \(\downarrow 29\%\) \\ \bottomrule
\end{tabular}
\end{table}
\FloatBarrier

\FloatBarrier
\subsection{Sensitivity analysis} 
\label{sec:sensitivity}
\begin{figure}
    \centering
    \includegraphics[width=\textwidth]{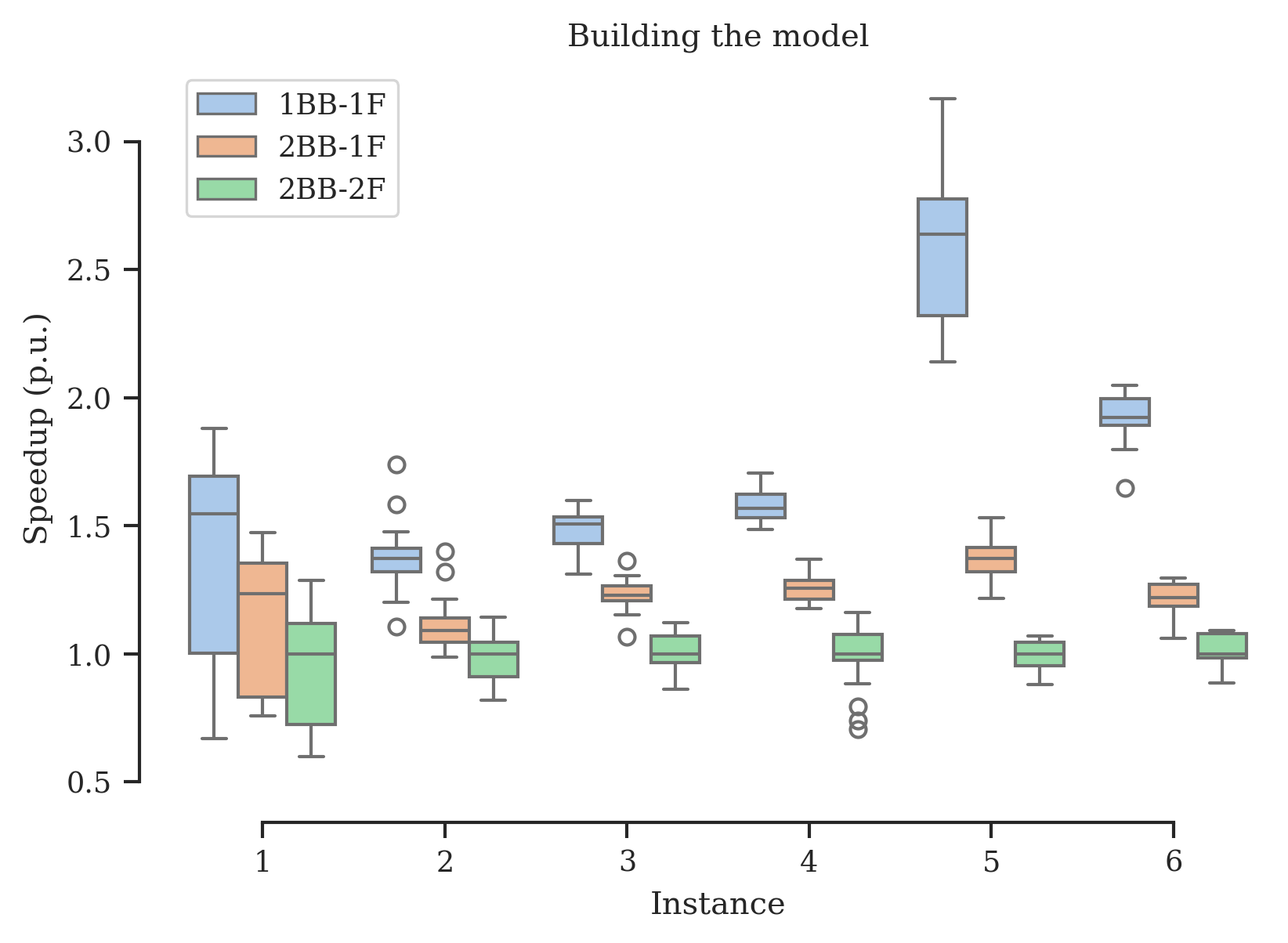}
    \caption{Speedups distribution comparison for the time to build}
    \label{fig:speedup-distribution-time-to-build}
\end{figure}

\begin{figure}
    \centering
    \includegraphics[width=\textwidth]{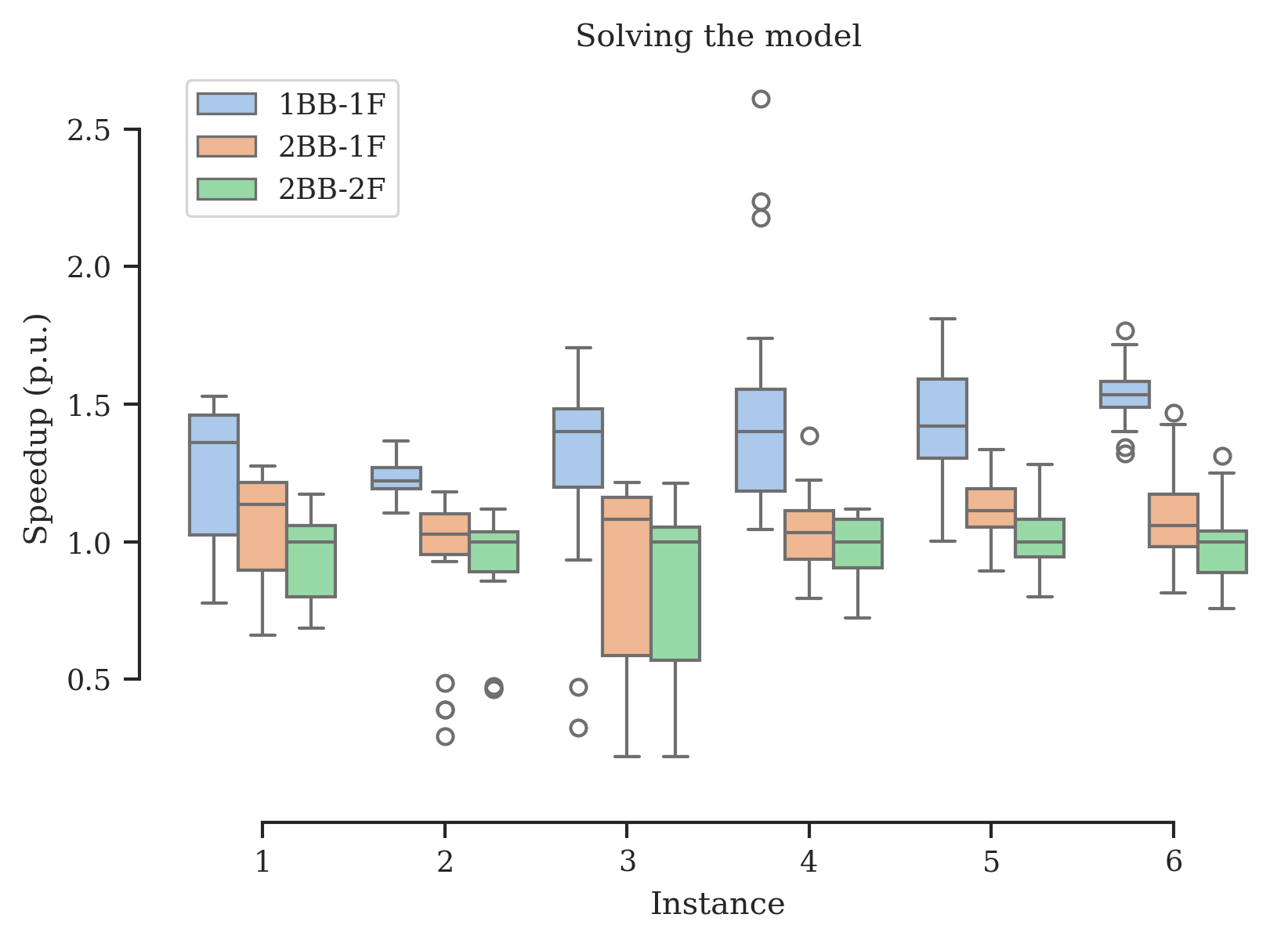}
    \caption{Speedups distribution comparison for the total time to solve}
    \label{fig:speedup-distribution-total-time-to-solve}
\end{figure}

It is important to comment that commercial optimisation solvers, such as Gurobi, can yield different results based on the "seed" parameter used \citep{Tejada2020}. Therefore, relying on a single run is not the best procedure for comparing different modelling approaches. Furthermore, the time it takes to build the model can vary due to random variations in CPU processing, which can also lead to differences between each run. To address this, we conducted a sensitivity analysis with 30 different seeds for each approach and optimisation instance, calculating median values for each group to determine if there are any statistically significant differences. We chose the median value as it is a central measure unaffected by data outliers. Figure \ref{fig:speedup-distribution-time-to-build} shows the distribution of speedups and their means with the 2BB-2F approach as a reference for all results in each instance. Similarly, Figure \ref{fig:speedup-distribution-total-time-to-solve} shows the results for solving time. Both figures suggest a difference between the 1BB-1F approach and both approaches using two building blocks, i.e., 2BB-2F and 2BB-1F. However, the 2BB-2F and 2BB-1F approaches are closer in the distribution of values. Section \ref{sec:statistical-analysis} provides a statistical analysis of the results to determine if there's enough statistical evidence to support these hypotheses.

\FloatBarrier
\subsection{Statistical analysis} 
\label{sec:statistical-analysis}

To determine if there is a significant difference between the modelling approaches, we use a two-sample t-test. The null hypothesis states that the average values of two related or repeated samples are identical. In contrast, the alternative hypothesis suggests that the underlying distributions of the samples have unequal means.

The variability in the results stems from the random number seed used by the solver, which leads the solver to take different solution paths. Therefore, we can assume for the t-test that the results follow a normal probability distribution, the variances of the results are equal, and each individual in the population has an equal probability of being selected in the sample.

We use 2BB-2F as the reference and compare it individually with  2BB-1F and 1BB-1F. Here, we assume a significance level of $\alpha = 0.05$. Table \ref{tab:t-test-nodal-2-flows} shows the test results for each instance and compared approach. If the p-value is smaller than the threshold, we reject the null hypothesis of equal averages. Then, there is a statistically significant difference between the mean times of 2BB-2F and the others (2BB-1F and 1BB-1F).

The results presented in Table \ref{tab:t-test-nodal-2-flows} indicate that 1BB-1F  has a significantly different mean value than 2BB-2F in all instances. Additionally, Figure \ref{fig:speedup-distribution-total-time-to-solve} shows that the average speedup values are higher for 1BB-1F in all cases. This implies that while 2BB-2F could potentially be faster than 1BB-1F in a particular instance, there is statistical evidence that, on average, 1BB-1F is faster. As for 2BB-1F, no statistical evidence suggests it is faster than 2BB-F for half of the instances. However, it is worth noting that the larger case studies in the test, instances 5 and 6, show statistical evidence that the 2BB-1F approach is faster, on average, than the 2BB-2F.

We are left with the question of how the 1BB-1F approach compares to the 2BB-1F approach. We also use a two-sample t-test with the 2BB-1F approach as the reference to compare these approaches. The results in Table \ref{tab:t-test-nodal-1-flow} indicate that the 1BB-1F approach is statistically faster for all instances, as shown in Figure \ref{fig:speedup-distribution-total-time-to-solve}.

\begin{table}[ht]
\centering
\caption{T-test results using 2BB-2F as reference}
\label{tab:t-test-nodal-2-flows}
\begin{tabular}{cccc}
\toprule
\textbf{Compared Approach} & \textbf{Instance} & \textbf{t-statistic} & \textbf{p-value} \\ \hline
       & 1  & 2.11  & 0.04 \\
       & 2  & -0.37 & \textbf{0.72} \\
2BB-1F & 3  & 0.34  & \textbf{0.74} \\
       & 4  & 1.74  & \textbf{0.09} \\
       & 5  & 4.47  & 0.00 \\ 
       & 6  & 2.60  & 0.01 \\ \midrule
       & 1  & 6.68  & 0.00 \\ 
       & 2  & 4.83  & 0.00 \\ 
1BB-1F & 3  & 4.78  & 0.00 \\ 
       & 4  & 9.49  & 0.00 \\ 
       & 5  & 11.48 & 0.00 \\ 
       & 6  & 15.08 & 0.00 \\ \bottomrule
\end{tabular}
\end{table}

\begin{table}[ht]
\centering
\caption{T-test results using 2BB-1F as reference}
\label{tab:t-test-nodal-1-flow}
\begin{tabular}{cccc}
\toprule
\textbf{Compared Approach} & \textbf{Instance} & \textbf{t-statistic} & \textbf{p-value} \\ \hline
        & 1 & 6.68  & 0.00 \\ 
        & 2 & 4.83  & 0.00 \\ 
1BB-1F  & 3 & 4.78  & 0.00 \\ 
        & 4 & 9.49  & 0.00 \\ 
        & 5 & 11.48 & 0.00 \\ 
        & 6 & 15.08 & 0.00 \\ \bottomrule
\end{tabular}
\end{table}

\FloatBarrier
\subsection{Discussion}
\label{sec:discussion}
Nowadays, energy system modellers have access to several energy models that can be used to develop case studies. \citet{LAVENEZIANA2023} present a critical review of energy planning models, comparing their capabilities and performance. However, this review does not explore the available modelling options described in this paper. Therefore, Table \ref{tab:models-overview} provides a general overview of the main available modelling approaches in a sample of these energy system models. Models focusing on power systems mainly consist of two building blocks, \textit{nodes} and energy \textit{assets}. In contrast, most recent models that focus on multi-sector analysis, like Calliope and SpineOpt, have three building blocks: \textit{nodes}, energy \textit{assets}, and \textit{connections}. The reason behind having more building blocks is to create a flexible model that can adapt to different energy asset configurations and represent different sectors. However, a flexible model can increase the number of variables and constraints, which can come with a computation cost. The results in the illustrative example in this paper estimate the impact in terms of speedup. As a general recommendation for energy system modellers, it is better to have fewer variables, especially for large-scale case studies.

In order to reduce the size of the model, instead of reformulating the model using a single building block as shown in this paper, some models provide the option to select specific methods tailored to the application. These methods also aim to reduce the number of variables and constraints. For instance, IRENA FlexTool gives users a selection of conversion and transfer methods, and some of those methods allow the use of fewer variables and constraints when the model is created, which also improves the solution process in the solver. FlexTool can represent the example problem in 1BB-1F format with one extra variable using available methods and in full 1BB-1F format using carefully formulated user constraints, which can be entered as data. 

Another example comes from PyPSA, a 2BB-1F approach model. In PyPSA, the bus serves as the fundamental building block, to which all energy assets—such as loads, generators, storage units, lines, transformers, and links—are connected. This design means that a direct connection between two energy assets is not possible, which distinguishes it from the 1BB-1F model. However, PyPSA provides the option to create custom constraints, allowing for more complex configurations, like those illustrated in the case study, while reducing the total number of constraints required.
In both cases, Flextool and PyPSA, the benefit comes from choosing a modelling approach, method, or custom constraint that allows fewer variables and constraints to represent the same energy system. The advantage and strength of the proposed 1BB-1F  in this paper is that it is a generic way to formulate the problem without depending on the definition of tailor-made methods or constraints that reduce model size depending on the user settings in a given case study.

The proposal for using only one building block with one flow variable (1BB-1F) is a new way of utilising the network graph structure of energy systems. While this approach has been found to offer several advantages, it also poses two primary challenges and limitations. First and foremost, it involves an extra layer of abstraction as the proposal can establish direct connections between energy assets, which may not be intuitive compared to traditional methods that use more building blocks like nodes and connections. However, this level of abstraction does not hinder the modellers from representing assets that perform the functions of nodes or connections, as evidenced in the case study where the same model, where we used TulipaEnergyModel to represent all the modelling approaches. Secondly, it is not always possible to obtain all the benefits of asset-to-asset connections (1BB-1F) in situations where constraints to model gas pressure, heat transfer, or DCOPF with losses for transmission or distribution networks are considered. Modelling these situations results in variations between an energy asset's incoming and outgoing flows, which partly reduces the potential savings on the model size of the 1BB-1F approach. In such cases, it becomes necessary to include two or four flow variables anyway to represent the situation accurately. Although the proposed 1BB-1F approach allows a direct connection between assets eliminating unnecessary elements in between, with their corresponding variables and constraints , it needs more variables to correctly represent more elaborated flows. Nevertheless, 1BB-1F  helps  large-scale optimisation problems that can be simplified.

\begin{table}[ht]
\centering
\caption{Overview of available modelling approaches in a sample of energy system models}
\label{tab:models-overview}
\begin{tabularx}{\textwidth}{c*{3}{>{\centering\arraybackslash}X}}
\toprule
\textbf{Model} & \textbf{Main focus} &  \textbf{Year} & \textbf{Approach} \\ 
\midrule 
Backbone \citep{Helistö2019}                    &  Multi-sector  & 2016 & 2BB-1F  \\	
Balmorel \citep{Wiese2018}                      &  Multi-sector	 & 2018 & 2BB-2F  \\	
Calliope \citep{Pfenninger2018}                 &  Multi-sector	 & 2018 & 3BB-4F  \\	
COMPETES-TNO \citep{ozdemir2019a}               &  Power Systems & 2004 & 2BB-2F  \\		
FlexTool \citep{FlexTool}                       &  Power Systems & 2021 & 2BB-1F \footnotemark{}  \\		
GenX \citep{Jenkins_GenX_2022}                  &  Power Systems & 2022 & 2BB-1F  \\		
OSEMOSYS \citep{OSeMOSYS}                       &  Multi-sector  & 2016 & 2BB-2F  \\		
Plexos \citep{Plexos}                           &  Multi-sector  & 2000 & 2BB-1F  \\		
PowerSystems \citep{LARA2021100747}             &  Power Systems & 2021 & 2BB-1F  \\		
PyPSA \citep{PyPSA}                             &  Multi-sector  & 2018 & 2BB-1F \footnotemark[\value{footnote}] \\ 
SpineOpt \citep{Ihlemann2022}                   &  Multi-sector  & 2022 & 3BB-4F  \\ 
TIMES \citep{TIMES_2016}                        &  Multi-sector  & 2016 & 2BB-2F  \\ 
TulipaEnergyModel \citep{TulipaEnergyModel2023} &  Multi-sector  & 2023 & 1BB-1F  \\ 
\bottomrule
\end{tabularx}
\end{table}

\footnotetext{Models with extra methods or custom constraints that allow for reducing the number of variables and constraints, depending on the modelling needs.}


\FloatBarrier
\section{Conclusion} 
\label{sec:conclusion}

This paper demonstrates that large-scale LP energy system models can be reformulated to reduce computational burden without compromising fidelity. We introduce a graph-based formulation using a single building block—the energy asset—and show that this approach leads to a more compact and computationally efficient LP formulation. Compared with four traditional approaches that rely on multiple building blocks and connections, our method reduces variables and constraints by 26\% and 35\%, respectively, and achieves an average speedup of 1.27× in solving time. These improvements scale with problem size, offering substantial benefits for large-scale, high-resolution studies.

Our findings challenge the widespread belief that LP models cannot be improved without a loss of accuracy. Even though solvers can eliminate some redundancy during presolve, we show that providing a compact, well-formulated model leads to faster building and solving times. This is particularly valuable for researchers and practitioners working on regional or multi-sector energy system models.

Beyond technical benefits, this work advocates for more conscious modelling practices. Flexible models should not only accommodate various configurations but do so efficiently. The proposed approach enables scalable, high-fidelity modelling using existing solvers (software) and computers (hardware), making it easier to explore future energy pathways with greater resolution and technical details. We hope this research encourages the energy modelling community to reconsider assumptions about LP formulations and adopt modelling strategies that balance flexibility, fidelity, and performance.

\section{CRediT author statement}
\textbf{Diego A. Tejada-Arango}: Methodology, Writing- Original draft preparation, Formal analysis, Software. \textbf{Germán Morales-España}: Methodology, Validation, Writing- Reviewing and Editing, Funding acquisition. \textbf{Juha Kiviluoma}: Conceptualization, Writing- Reviewing and Editing, Funding acquisition.

\section{Declaration of Generative AI and AI-assisted technologies in the writing process}
During the preparation of this work, the authors used Grammarly to improve readability and language. After using this tool/service, the authors reviewed and edited the content as needed and took full responsibility for the publication's content.

\section{Acknowledgments}

This research received funding from the European Climate, Infrastructure and Environment Executive Agency under the European Union’s HORIZON Research and Innovation Actions under grant agreement N°101095998. In addition, the Dutch Research Council (NWO) also partially funded this research under grant number ESI.2019.008.

\textbf{Disclaimer}: Views and opinions expressed are those of the author(s) only and those of the European Union or NWO. Neither the European Union nor the granting authority can be held responsible.

We would like to express our gratitude to the three anonymous reviewers for their valuable feedback, which significantly enhanced the initial version of this manuscript. Additionally, we extend our thanks to Fabian Neumann for his assistance in better understanding PyPSA's formulation and capabilities, which were essential for the comments on PyPSA in this paper.

\section{Supplementary material}
\label{sec:supplementary_material}
All the code to run the experiments and the raw results in this paper are available at the following link:

\url{https://github.com/datejada/experiments-flexible-connection}

The formulation presented in this paper is a simplified version of the TulipaEnergyModel.jl. The complete formulation, including all features of the model, can be found at the following link:

\url{https://tulipaenergy.github.io/TulipaEnergyModel.jl/stable/}

\section*{Annex A. Extended formulation}
\label{sec:annex-a}
Section \ref{sec:generic-formulation} shows a basic formulation for the proposed approach. In this section, we illustrate how to include other modelling features; for instance, here we show the equations to model the DC power flow and unit commitment problem.

\subsection*{A.1. DC power flow}

For a flow connecting assets $a^{\text{from}}$ and $a^{\text{to}}$, which belongs to the set $\mathcal{F}^{\text{dc-opf}}$, the power flow constraints are as follows:

\begin{equation}
\begin{aligned}
f_{(a^{\text{from}},a^{\text{to}}),t} = S^{\text{base}} \cdot \frac{\theta_{a^{\text{from}},t} - \theta_{a^{\text{to}},t}}{X_{(a^{\text{from}},a^{\text{to}})}} \quad \forall (a^{\text{from}},a^{\text{to}}) \in \mathcal{F}^{\text{dc-opf}}, \forall t \in \mathcal{T}
\end{aligned}
\end{equation}

Where:

\begin{itemize}
    \item[] $\mathcal{F}^{\text{dc-opf}} \subseteq \mathcal{F}$: subset of flows that have DC power flow constraints
    \item[] $\theta_{a,t}$: voltage angle of asset $a$ in timestep $t$. Notice that this asset is playing a role of a bus in the traditional energy modelling approach.
    \item[] $X_{(a^{\text{from}},a^{\text{to}})}$: per unit reactance between two assets, $a^{\text{from}}$ and $a^{\text{to}}$
    \item[] $S^{base}$: base power, e.g., 100 MVA 
\end{itemize}

\subsection*{A.2. Unit commitment}
Assets within the subset $\mathcal{A}^{\text{uc}}$ will contain the unit commitment constraints in the model. The following constraints show a basic formulation using the unit commitment variable $u$ for the units on. However, more complex constraints with startup and shutdown variables, e.g., \citep{gentile_tight_2016}, can also be added to the model.

\begin{equation}
\label{eq:min_oper_point}
\hat{f}_{a,t} = \sum_{(a,a^{\text{to}}) \in \mathcal{F}} f_{(a,a^{\text{to}}),t} - \underline{P}_a \cdot u_{a,t}  \quad \forall a \in \mathcal{A}^{\text{uc}}, \forall t \in \mathcal{T}
\end{equation}

\begin{equation}
\label{eq:uc_limit}
u_{a,t} \leq U^0_{a} + i_{a}  \quad \forall a \in \mathcal{A}^{\text{uc}}, \forall t \in \mathcal{T}
\end{equation}

\begin{equation}
\label{eq:max_flow_above_min}
\hat{f}_{a,t} \leq \left(\overline{P}_a - \underline{P}_a \right) \cdot u_{a,t}  \quad \forall a \in \mathcal{A}^{\text{uc}}, \forall t \in \mathcal{T}
\end{equation}

\begin{equation}
\label{eq:min_flow_above_min}
\hat{f}_{a,t} \geq 0  \quad \forall a \in \mathcal{A}^{\text{uc}}, \forall t \in \mathcal{T}
\end{equation}

Where:

\begin{itemize}
    \item[] $\mathcal{A}^{\text{uc}}  \subseteq \mathcal{A}$: subset of assets with unit commitment 
    \item[] $\hat{f}_{a,t}$: flow above the minimum operating point of asset $a$ in timestep $t$
    \item[] $f_{(a,a^{\text{to}}),t}$: flow variable from asset $a$ to asset $a^{\text{to}}$ in timestep $t$    
    \item[] $u_{a,t}$: units on variable of $a$ in timestep $t$
    \item[] $i_{a}$: investment unit variable of asset $a$    
    \item[] $U^0_{a}$: initial number of units parameter of the asset $a$    
    \item[] $\overline{P}_{a}$: maximum capacity parameter of the asset $a$
    \item[] $\underline{P}_{a}$: minimum capacity parameter of the asset $a$       
\end{itemize}

Equation (\ref{eq:min_oper_point}) defines the flow above the asset's minimum operating point, and  (\ref{eq:uc_limit}) sets the limits on the variable units. Equation (\ref{eq:max_flow_above_min}) specifies the maximum output flow above the minimum operating point, while  (\ref{eq:min_flow_above_min}) indicates the minimum output flow above the minimum operating point.

\bibliographystyle{elsarticle-num-names} 
\bibliography{cas-refs}

\begin{thebibliography}{27}
\expandafter\ifx\csname natexlab\endcsname\relax\def\natexlab#1{#1}\fi
\providecommand{\url}[1]{\texttt{#1}}
\providecommand{\href}[2]{#2}
\providecommand{\path}[1]{#1}
\providecommand{\DOIprefix}{doi:}
\providecommand{\ArXivprefix}{arXiv:}
\providecommand{\URLprefix}{URL: }
\providecommand{\Pubmedprefix}{pmid:}
\providecommand{\doi}[1]{\href{http://dx.doi.org/#1}{\path{#1}}}
\providecommand{\Pubmed}[1]{\href{pmid:#1}{\path{#1}}}
\providecommand{\bibinfo}[2]{#2}
\ifx\xfnm\relax \def\xfnm[#1]{\unskip,\space#1}\fi
\bibitem[{Kuhn(2010)}]{Kuhn2010}
\bibinfo{author}{H.~W. Kuhn}, \bibinfo{title}{The Hungarian Method for the Assignment Problem}, \bibinfo{publisher}{Springer Berlin Heidelberg}, \bibinfo{address}{Berlin, Heidelberg}, \bibinfo{year}{2010}, pp. \bibinfo{pages}{29--47}.
\bibitem[{Cherkassky and Goldberg(1997)}]{Cherkassky1997}
\bibinfo{author}{B.~V. Cherkassky}, \bibinfo{author}{A.~V. Goldberg},
\newblock \bibinfo{title}{On implementing the push—relabel method for the maximum flow problem},
\newblock \bibinfo{journal}{Algorithmica} \bibinfo{volume}{19} (\bibinfo{year}{1997}) \bibinfo{pages}{390--410}. \URLprefix \url{https://doi.org/10.1007/PL00009180}. \DOIprefix\doi{10.1007/PL00009180}.
\bibitem[{TNO(2020)}]{ESDL2020}
\bibinfo{author}{TNO}, \bibinfo{title}{Energy system description language}, \bibinfo{year}{2020}. \URLprefix \url{https://energytransition.gitbook.io/esdl/esdl-concepts/design-principles}, \bibinfo{note}{accessed on March 19th, 2024}.
\bibitem[{Wiese et~al.(2018)Wiese, Bramstoft, Koduvere, {Pizarro Alonso}, Balyk, Kirkerud, Åsa Grytli~Tveten, Bolkesjø, Münster, and Ravn}]{Wiese2018}
\bibinfo{author}{F.~Wiese}, \bibinfo{author}{R.~Bramstoft}, \bibinfo{author}{H.~Koduvere}, \bibinfo{author}{A.~{Pizarro Alonso}}, \bibinfo{author}{O.~Balyk}, \bibinfo{author}{J.~G. Kirkerud}, \bibinfo{author}{Åsa Grytli~Tveten}, \bibinfo{author}{T.~F. Bolkesjø}, \bibinfo{author}{M.~Münster}, \bibinfo{author}{H.~Ravn},
\newblock \bibinfo{title}{Balmorel open source energy system model},
\newblock \bibinfo{journal}{Energy Strategy Reviews} \bibinfo{volume}{20} (\bibinfo{year}{2018}) \bibinfo{pages}{26--34}. \URLprefix \url{https://www.sciencedirect.com/science/article/pii/S2211467X18300038}. \DOIprefix\doi{https://doi.org/10.1016/j.esr.2018.01.003}.
\bibitem[{{Energy Exemplar}(2000)}]{Plexos}
\bibinfo{author}{{Energy Exemplar}}, \bibinfo{title}{Plexos}, \bibinfo{year}{2000}. \URLprefix \url{https://www.energyexemplar.com/plexos}, \bibinfo{note}{accessed on March 11th, 2024}.
\bibitem[{Morales-España et~al.(2024)Morales-España, Hernández-Serna, Tejada-Arango, and Weeda}]{MoralesEspaña2024}
\bibinfo{author}{G.~Morales-España}, \bibinfo{author}{R.~Hernández-Serna}, \bibinfo{author}{D.~A. Tejada-Arango}, \bibinfo{author}{M.~Weeda},
\newblock \bibinfo{title}{Impact of large-scale hydrogen electrification and retrofitting of natural gas infrastructure on the european power system},
\newblock \bibinfo{journal}{International Journal of Electrical Power \& Energy Systems} \bibinfo{volume}{155} (\bibinfo{year}{2024}) \bibinfo{pages}{109686}. \URLprefix \url{https://www.sciencedirect.com/science/article/pii/S0142061523007433}. \DOIprefix\doi{https://doi.org/10.1016/j.ijepes.2023.109686}.
\bibitem[{Gea-Bermúdez et~al.(2021)Gea-Bermúdez, Jensen, Münster, Koivisto, Kirkerud, kuang Chen, and Ravn}]{Gea-Bermudez2021}
\bibinfo{author}{J.~Gea-Bermúdez}, \bibinfo{author}{I.~G. Jensen}, \bibinfo{author}{M.~Münster}, \bibinfo{author}{M.~Koivisto}, \bibinfo{author}{J.~G. Kirkerud}, \bibinfo{author}{Y.~kuang Chen}, \bibinfo{author}{H.~Ravn},
\newblock \bibinfo{title}{The role of sector coupling in the green transition: A least-cost energy system development in northern-central europe towards 2050},
\newblock \bibinfo{journal}{Applied Energy} \bibinfo{volume}{289} (\bibinfo{year}{2021}) \bibinfo{pages}{116685}. \URLprefix \url{https://www.sciencedirect.com/science/article/pii/S0306261921002130}. \DOIprefix\doi{https://doi.org/10.1016/j.apenergy.2021.116685}.
\bibitem[{Kiviluoma et~al.(2022)Kiviluoma, Tupala, Soininen, and {International Renewable Energy Agency (IRENA)}}]{FlexTool}
\bibinfo{author}{J.~Kiviluoma}, \bibinfo{author}{A.~Tupala}, \bibinfo{author}{A.~Soininen}, \bibinfo{author}{{International Renewable Energy Agency (IRENA)}}, \bibinfo{title}{{IRENA FlexTool}}, \bibinfo{year}{2022}. \URLprefix \url{https://github.com/irena-flextool/flextool}.
\bibitem[{Ihlemann et~al.(2022)Ihlemann, Kouveliotis-Lysikatos, Huang, Dillon, O’Dwyer, Rasku, Marin, Poncelet, and Kiviluoma}]{Ihlemann2022}
\bibinfo{author}{M.~Ihlemann}, \bibinfo{author}{I.~Kouveliotis-Lysikatos}, \bibinfo{author}{J.~Huang}, \bibinfo{author}{J.~Dillon}, \bibinfo{author}{C.~O’Dwyer}, \bibinfo{author}{T.~Rasku}, \bibinfo{author}{M.~Marin}, \bibinfo{author}{K.~Poncelet}, \bibinfo{author}{J.~Kiviluoma},
\newblock \bibinfo{title}{Spineopt: A flexible open-source energy system modelling framework},
\newblock \bibinfo{journal}{Energy Strategy Reviews} \bibinfo{volume}{43} (\bibinfo{year}{2022}) \bibinfo{pages}{100902}. \URLprefix \url{https://www.sciencedirect.com/science/article/pii/S2211467X22000955}. \DOIprefix\doi{https://doi.org/10.1016/j.esr.2022.100902}.
\bibitem[{Markensteijn et~al.(2020)Markensteijn, Romate, and Vuik}]{Markensteijn2020}
\bibinfo{author}{A.~Markensteijn}, \bibinfo{author}{J.~Romate}, \bibinfo{author}{C.~Vuik},
\newblock \bibinfo{title}{A graph-based model framework for steady-state load flow problems of general multi-carrier energy systems},
\newblock \bibinfo{journal}{Applied Energy} \bibinfo{volume}{280} (\bibinfo{year}{2020}) \bibinfo{pages}{115286}. \URLprefix \url{https://www.sciencedirect.com/science/article/pii/S0306261920307984}. \DOIprefix\doi{https://doi.org/10.1016/j.apenergy.2020.115286}.
\bibitem[{Klotz and Newman(2013)}]{KLOTZ20131}
\bibinfo{author}{E.~Klotz}, \bibinfo{author}{A.~M. Newman},
\newblock \bibinfo{title}{Practical guidelines for solving difficult linear programs},
\newblock \bibinfo{journal}{Surveys in Operations Research and Management Science} \bibinfo{volume}{18} (\bibinfo{year}{2013}) \bibinfo{pages}{1--17}. \URLprefix \url{https://www.sciencedirect.com/science/article/pii/S1876735412000189}. \DOIprefix\doi{https://doi.org/10.1016/j.sorms.2012.11.001}.
\bibitem[{Lustig et~al.(1991)Lustig, Mulvey, and Carpenter}]{Irvin_1991}
\bibinfo{author}{I.~J. Lustig}, \bibinfo{author}{J.~M. Mulvey}, \bibinfo{author}{T.~J. Carpenter},
\newblock \bibinfo{title}{Formulating two-stage stochastic programs for interior point methods},
\newblock \bibinfo{journal}{Operations Research} \bibinfo{volume}{39} (\bibinfo{year}{1991}) \bibinfo{pages}{757--770}. \URLprefix \url{https://doi.org/10.1287/opre.39.5.757}. \DOIprefix\doi{10.1287/opre.39.5.757}. \href{http://arxiv.org/abs/https://doi.org/10.1287/opre.39.5.757}{{\tt arXiv:https://doi.org/10.1287/opre.39.5.757}}.
\bibitem[{Brown et~al.(2018)Brown, H\"orsch, and Schlachtberger}]{PyPSA}
\bibinfo{author}{T.~Brown}, \bibinfo{author}{J.~H\"orsch}, \bibinfo{author}{D.~Schlachtberger},
\newblock \bibinfo{title}{{PyPSA: Python for Power System Analysis}},
\newblock \bibinfo{journal}{Journal of Open Research Software} \bibinfo{volume}{6} (\bibinfo{year}{2018}). \URLprefix \url{https://doi.org/10.5334/jors.188}. \DOIprefix\doi{10.5334/jors.188}. \href{http://arxiv.org/abs/1707.09913}{{\tt arXiv:1707.09913}}.
\bibitem[{Howells et~al.(2011)Howells, Rogner, Strachan, Heaps, Huntington, Kypreos, Hughes, Silveira, DeCarolis, Bazillian, and Roehrl}]{OSeMOSYS}
\bibinfo{author}{M.~Howells}, \bibinfo{author}{H.~Rogner}, \bibinfo{author}{N.~Strachan}, \bibinfo{author}{C.~Heaps}, \bibinfo{author}{H.~Huntington}, \bibinfo{author}{S.~Kypreos}, \bibinfo{author}{A.~Hughes}, \bibinfo{author}{S.~Silveira}, \bibinfo{author}{J.~DeCarolis}, \bibinfo{author}{M.~Bazillian}, \bibinfo{author}{A.~Roehrl},
\newblock \bibinfo{title}{Osemosys: The open source energy modeling system: An introduction to its ethos, structure and development},
\newblock \bibinfo{journal}{Energy Policy} \bibinfo{volume}{39} (\bibinfo{year}{2011}) \bibinfo{pages}{5850--5870}. \URLprefix \url{https://www.sciencedirect.com/science/article/pii/S0301421511004897}. \DOIprefix\doi{https://doi.org/10.1016/j.enpol.2011.06.033}, \bibinfo{note}{sustainability of biofuels}.
\bibitem[{Tejada-Arango et~al.(2023)Tejada-Arango, Morales-España, Clisby, Wang, Soares~Siqueira, Ali, Soucasse, and Neustroev}]{TulipaEnergyModel2023}
\bibinfo{author}{D.~A. Tejada-Arango}, \bibinfo{author}{G.~Morales-España}, \bibinfo{author}{L.~Clisby}, \bibinfo{author}{N.~Wang}, \bibinfo{author}{A.~Soares~Siqueira}, \bibinfo{author}{S.~Ali}, \bibinfo{author}{L.~Soucasse}, \bibinfo{author}{G.~Neustroev}, \bibinfo{title}{{Tulipa Energy Model}}, \bibinfo{year}{2023}. \URLprefix \url{https://github.com/TulipaEnergy/TulipaEnergyModel.jl}.
\bibitem[{Bezanson et~al.(2017)Bezanson, Edelman, Karpinski, and Shah}]{Bezanson_Julia_A_fresh_2017}
\bibinfo{author}{J.~Bezanson}, \bibinfo{author}{A.~Edelman}, \bibinfo{author}{S.~Karpinski}, \bibinfo{author}{V.~B. Shah},
\newblock \bibinfo{title}{{Julia: A fresh approach to numerical computing}},
\newblock \bibinfo{journal}{SIAM Review} \bibinfo{volume}{59} (\bibinfo{year}{2017}) \bibinfo{pages}{65--98}. \DOIprefix\doi{10.1137/141000671}.
\bibitem[{Lubin et~al.(2023)Lubin, Dowson, {Dias Garcia}, Huchette, Legat, and Vielma}]{Lubin2023}
\bibinfo{author}{M.~Lubin}, \bibinfo{author}{O.~Dowson}, \bibinfo{author}{J.~{Dias Garcia}}, \bibinfo{author}{J.~Huchette}, \bibinfo{author}{B.~Legat}, \bibinfo{author}{J.~P. Vielma},
\newblock \bibinfo{title}{{JuMP} 1.0: {R}ecent improvements to a modeling language for mathematical optimization},
\newblock \bibinfo{journal}{Mathematical Programming Computation}  (\bibinfo{year}{2023}). \DOIprefix\doi{10.1007/s12532-023-00239-3}.
\bibitem[{Fairbanks et~al.(2021)Fairbanks, Besan{\c{c}}on, Simon, Hoffiman, Eubank, and Karpinski}]{Graphs2021}
\bibinfo{author}{J.~Fairbanks}, \bibinfo{author}{M.~Besan{\c{c}}on}, \bibinfo{author}{S.~Simon}, \bibinfo{author}{J.~Hoffiman}, \bibinfo{author}{N.~Eubank}, \bibinfo{author}{S.~Karpinski}, \bibinfo{title}{Juliagraphs/graphs.jl: an optimized graphs package for the julia programming language}, \bibinfo{year}{2021}. \URLprefix \url{https://github.com/JuliaGraphs/Graphs.jl/}.
\bibitem[{Tejada-Arango et~al.(2020)Tejada-Arango, Lumbreras, Sánchez-Martín, and Ramos}]{Tejada2020}
\bibinfo{author}{D.~A. Tejada-Arango}, \bibinfo{author}{S.~Lumbreras}, \bibinfo{author}{P.~Sánchez-Martín}, \bibinfo{author}{A.~Ramos},
\newblock \bibinfo{title}{Which unit-commitment formulation is best? a comparison framework},
\newblock \bibinfo{journal}{IEEE Transactions on Power Systems} \bibinfo{volume}{35} (\bibinfo{year}{2020}) \bibinfo{pages}{2926--2936}. \DOIprefix\doi{10.1109/TPWRS.2019.2962024}.
\bibitem[{Laveneziana et~al.(2023)Laveneziana, Prussi, and Chiaramonti}]{LAVENEZIANA2023}
\bibinfo{author}{L.~Laveneziana}, \bibinfo{author}{M.~Prussi}, \bibinfo{author}{D.~Chiaramonti},
\newblock \bibinfo{title}{Critical review of energy planning models for the sustainable development at company level},
\newblock \bibinfo{journal}{Energy Strategy Reviews} \bibinfo{volume}{49} (\bibinfo{year}{2023}) \bibinfo{pages}{101136}. \URLprefix \url{https://www.sciencedirect.com/science/article/pii/S2211467X2300086X}. \DOIprefix\doi{https://doi.org/10.1016/j.esr.2023.101136}.
\bibitem[{Helistö et~al.(2019)Helistö, Kiviluoma, Ikäheimo, Rasku, Rinne, O’Dwyer, Li, and Flynn}]{Helistö2019}
\bibinfo{author}{N.~Helistö}, \bibinfo{author}{J.~Kiviluoma}, \bibinfo{author}{J.~Ikäheimo}, \bibinfo{author}{T.~Rasku}, \bibinfo{author}{E.~Rinne}, \bibinfo{author}{C.~O’Dwyer}, \bibinfo{author}{R.~Li}, \bibinfo{author}{D.~Flynn},
\newblock \bibinfo{title}{Backbone—an adaptable energy systems modelling framework},
\newblock \bibinfo{journal}{Energies} \bibinfo{volume}{12} (\bibinfo{year}{2019}). \URLprefix \url{https://www.mdpi.com/1996-1073/12/17/3388}. \DOIprefix\doi{10.3390/en12173388}.
\bibitem[{Pfenninger and Pickering(2018)}]{Pfenninger2018}
\bibinfo{author}{S.~Pfenninger}, \bibinfo{author}{B.~Pickering},
\newblock \bibinfo{title}{Calliope: a multi-scale energy systems modelling framework},
\newblock \bibinfo{journal}{Journal of Open Source Software} \bibinfo{volume}{3} (\bibinfo{year}{2018}) \bibinfo{pages}{825}. \URLprefix \url{https://doi.org/10.21105/joss.00825}. \DOIprefix\doi{10.21105/joss.00825}.
\bibitem[{Özge Özdemir et~al.(2019)Özge Özdemir, Hobbs, van Hout, and Koutstaal}]{ozdemir2019a}
\bibinfo{author}{Özge Özdemir}, \bibinfo{author}{B.~F. Hobbs}, \bibinfo{author}{M.~van Hout}, \bibinfo{author}{P.~Koutstaal},
\newblock \bibinfo{title}{Capacity vs energy subsidies for renewables: Benefits and costs for the 2030 eu power market},
\newblock \bibinfo{journal}{Energy Policy Research Group, University of Cambridge}  (\bibinfo{year}{2019}). \URLprefix \url{http://www.jstor.org/stable/resrep30402}.
\bibitem[{Jenkins et~al.(2022)Jenkins, Sepulveda, Mallapragada, Patankar, Schwartz, Schwartz, Chakrabarti, Xu, Morris, and Sepulveda}]{Jenkins_GenX_2022}
\bibinfo{author}{J.~Jenkins}, \bibinfo{author}{N.~Sepulveda}, \bibinfo{author}{D.~Mallapragada}, \bibinfo{author}{N.~Patankar}, \bibinfo{author}{A.~Schwartz}, \bibinfo{author}{J.~Schwartz}, \bibinfo{author}{S.~Chakrabarti}, \bibinfo{author}{Q.~Xu}, \bibinfo{author}{J.~Morris}, \bibinfo{author}{N.~Sepulveda}, \bibinfo{title}{{GenX}}, \bibinfo{year}{2022}. \URLprefix \url{https://github.com/GenXProject/GenX}. \DOIprefix\doi{10.5281/zenodo.6229425}.
\bibitem[{Lara et~al.(2021)Lara, Barrows, Thom, Krishnamurthy, and Callaway}]{LARA2021100747}
\bibinfo{author}{J.~D. Lara}, \bibinfo{author}{C.~Barrows}, \bibinfo{author}{D.~Thom}, \bibinfo{author}{D.~Krishnamurthy}, \bibinfo{author}{D.~Callaway},
\newblock \bibinfo{title}{Powersystems.jl — a power system data management package for large scale modeling},
\newblock \bibinfo{journal}{SoftwareX} \bibinfo{volume}{15} (\bibinfo{year}{2021}) \bibinfo{pages}{100747}. \URLprefix \url{https://www.sciencedirect.com/science/article/pii/S2352711021000765}. \DOIprefix\doi{https://doi.org/10.1016/j.softx.2021.100747}.
\bibitem[{Loulou et~al.(2016)Loulou, Goldstein, Kanudia, Lettila, Remme, Wright, Giannakidis, and Noble}]{TIMES_2016}
\bibinfo{author}{R.~Loulou}, \bibinfo{author}{G.~Goldstein}, \bibinfo{author}{A.~Kanudia}, \bibinfo{author}{A.~Lettila}, \bibinfo{author}{U.~Remme}, \bibinfo{author}{E.~Wright}, \bibinfo{author}{G.~Giannakidis}, \bibinfo{author}{K.~Noble}, \bibinfo{title}{{Documentation for the TIMES Model - Part I}}, \bibinfo{year}{2016}. \URLprefix \url{https://iea-etsap.org/index.php/documentation}.
\bibitem[{Gentile et~al.(2016)Gentile, Morales-España, and Ramos}]{gentile_tight_2016}
\bibinfo{author}{C.~Gentile}, \bibinfo{author}{G.~Morales-España}, \bibinfo{author}{A.~Ramos},
\newblock \bibinfo{title}{A tight {MIP} formulation of the unit commitment problem with start-up and shut-down constraints},
\newblock \bibinfo{journal}{EURO Journal on Computational Optimization} \bibinfo{volume}{5} (\bibinfo{year}{2016}) \bibinfo{pages}{177--201}. \URLprefix \url{https://link.springer.com/article/10.1007/s13675-016-0066-y}. \DOIprefix\doi{10.1007/s13675-016-0066-y}.

\end{thebibliography}





\end{document}